\documentclass[a4paper]{article}
\usepackage{graphicx}
\usepackage{hyperref}
\usepackage{latexsym}
\newtheorem{theorem}{Theorem}
\newtheorem{corollary}{Corollary}
\newtheorem{lemma}{Lemma}
\newtheorem{definition}{Definition}
\newtheorem{thesis}{Thesis}

\newtheorem{conclusion}{Conclusion}
\begin{document}
\title{\LARGE{Does resolving PvNP require a paradigm shift?} \\ \vspace{+1ex} \large{An investigation into the philosophical and mathematical significance of Aristotle's particularisation in the foundations of mathematics, logic and computability}}
\author{\normalsize{Bhupinder Singh Anand}}
\date{\small \textit {Draft of \today.}}
\maketitle

\begin{abstract}
I shall argue that a resolution of the PvNP problem requires building an iff bridge between the domain of provability and that of computability. The former concerns how a human intelligence decides the truth of number-theoretic relations, and is formalised by the first-order Peano Arithmetic PA following Dededekind's axiomatisation of Peano's Postulates. The latter concerns how a human intelligence computes the values of number-theoretic functions, and is formalised by the operations of a Turing Machine following Turing's analysis of computable functions. I shall show that such a bridge requires objective definitions of both an `algorithmic' interpretation of PA, and an `instantiational' interpretation of PA. I shall show that both interpretations are implicit in the definition of the subjectively defined `standard' interpretation of PA. However the existence of, and distinction between, the two objectively definable interpretations---and the fact that the former is sound whilst the latter is not---is obscured by the extraneous presumption under the `standard' interpretation of PA that Aristotle's particularisation must hold over the structure $\mathcal{N}$ of the natural numbers. I shall argue that recognising the falseness of this belief awaits a paradigm shift in our perception of the application of Tarski's analysis (of the concept of truth in the languages of the deductive sciences) to the `standard' interpretation of PA. I shall then show that an arithmetical formula $[F]$ is PA-provable if, and only if, $[F]$ interprets as true under an algorithmic interpretation of PA. I shall finally show how it then follows from G\"{o}del's construction of a formally `undecidable' arithmetical proposition that there is a Halting-type PA formula which---by Tarski's definitions---is algorithmically verifiable as true, but not algorithmically computable as true, under a sound interpretation of PA.
\end{abstract}

\section{Introduction}
\label{sec:2.1}
I define what it means for a number-theoretic function to be:

\begin{quote}
(i) Instantiationally computable;

\vspace{+.5ex}
(ii) Algorithmically computable.
\end{quote}

I argue that P$\neq$NP if a number-theoretic function is instantiationally computable but not algorithmically computable.

I then show in Lemma \ref{sec:2.3.lem.1} below that if Aristotle's particularisation is presumed valid over the structure $\mathcal{N}$ of the natural numbers---as is the case under the standard interpretation of PA---then it follows from the instantiational nature of the constructive definition of the G\"{o}del $\beta$-function\footnote{Introduced by Kurt G\"{o}del in Theorem VII (\cite{Go31}, pp.30-31) of his seminal 1931 paper on formally undecidable arithmetical propositions.} that a primitive recursive relation can be instantiationally equivalent to an arithmetical relation where the former is algorithmically computable as always true over $\mathcal{N}$ whilst the latter is instantiationally computable, but not algorithmically computable, as always true over $\mathcal{N}$.

\begin{quote}
\footnotesize{
I note that this argument cannot be formalised in ZF since functions are defined extensionally as mappings. Hence ZF cannot recognise that a primitive recursive relation may be instantiationally equivalent to, but computationally different from, an arithmetical relation where the former is algorithmically computable as always true over $\mathcal{N}$ whilst the latter is instantiationally computable, but not algorithmically computable, as always true over $\mathcal{N}$.
}
\end{quote}

I therefore conclude in Theorem \ref{sec:2.3.thm.1} that, if the standard interpretation of PA is presumed sound, then P$\neq$NP.

I then consider an algorithmic interpretation of PA that yields P$\neq$NP in Theorem \ref{sec:6.3.thm.2}, but which does not appeal to Aristotle's particularisation and which is provably sound.

\section{The PvNP problem}
\label{sec:1}

In a 2009 survey of the status of the P versus NP problem, Lance Fortnow wrote\footnote{\cite{Fo09}.}:

\begin{quote}
``\ldots in the mid-1980's, many believed that the quickly developing area of circuit complexity would soon settle the P versus NP problem, whether every algorithmic problem with efficiently verifiable solutions have efficiently computable solutions. But circuit complexity and other approaches to the problem have stalled and we have little reason to believe we will see a proof separating P from NP in the near future.

\ldots As we solve larger and more complex problems with greater computational power and cleverer algorithms, the problems we cannot tackle begin to stand out. The theory of NP-completeness helps us understand these limitations and the P versus NP problems begins to loom large not just as an interesting theoretical question in computer science, but as a basic principle that permeates all the sciences.

\ldots None of us truly understand the P versus NP problem, we have only begun to peel the layers around this increasingly complex question."
\end{quote}

\subsection{Equivalent definitions of P, NP and P$\neq$NP}
\label{sec:1.01}

The formal definition of the class P by Stephen Cook\footnote{\cite{Cook}.} admits a number-theoretic function $F$---viewed set-theoretically as defining (and defined by) a unique subset $L$ of the set $\Sigma^{*}$ of finite strings over some non-empty finite alphabet set $\Sigma$---in P if, and only if, some deterministic Turing machine TM accepts $L$ and runs in polynomial time.

\begin{quote}
\footnotesize{
In this investigation I interpret number-theoretic functions and relations over an infinite domain $\mathcal{D}$ as pre-Cantorian computational instructions (which may, or may not, be uniform) that, for any given sequence of allowable values to the variables in the function/relation, determine how the function/relation is to be evaluated---and whether, or not, the result of such evaluation yields a value (or values)---in the domain $\mathcal{D}$. I do not assume---as in Cantorian set theories---that the evaluations always determine a completed infinity (set) that can be referred to as a unique mathematical constant that identifies the function/relation in a mathematical language (or its interpretation) outside of the set theory in which the function/relation is defined.
}
\end{quote}

Fortnow describes the PvNP problem informally as follows:

\begin{quote}
``In 1965, Jack Edmonds \ldots suggested a formal definition of ``efficient computation" (runs in time a fixed polynomial of the input size). The class of problems with efficient solutions would later become known as P for ``Polynomial Time".

\ldots But many related problems do not seem to have such an efficient algorithm.

\ldots The collection of problems that have efficiently verifiable solutions is known as NP (for ``Nondeterministic Polynomial-Time" \dots).

So P=NP means that for every problem that has an efficiently verifiable solution, we can find that solution efficiently as well.

\ldots If a formula $\phi$ is not a tautology, we can give an easy proof of that fact by exhibiting an assignment of the variables that makes $\phi$ false. But if \ldots there are no short proofs of tautology that would imply P$\neq$NP."
\end{quote}

In an earlier paper presented to ICM 2002, Ran Raz explains\footnote{\cite{Ra02}.}:

\begin{quote}
``A Boolean formula $f(x_{1}, \ldots, x_{n})$ is a tautology if $f(x_{1}, \ldots, x_{n}) = 1$ for every $x_{1}, \ldots, x_{n}$. A Boolean formula $f(x_{1}, \ldots, x_{n})$ is unsatisfiable if $f(x_{1}, \ldots, x_{n}) = 0$ for every $x_{1}, \dots, x_{n}$. Obviously, $f$ is a tautology if and only if $\neg f$ is unsatisfiable.

Given a formula $f(x_{1}, \ldots, x_{n})$, one can decide whether or not $f$ is a tautology by checking all the possibilities for assignments to $x_{1}, \ldots, x_{n}$. However, the time needed for this procedure is exponential in the number of variables, and hence may be exponential in the length of the formula $f$.

\ldots P$\neq$NP is the central open problem in complexity theory and one of the most important open problems in mathematics today. The problem has thousands of equivalent formulations. One of these formulations is the following:

Is there a polynomial time algorithm $\mathcal{A}$ that gets as input a Boolean formula $f$ and outputs 1 if and only if $f$ is a tautology?

P$\neq$NP states that there is no such algorithm."
\end{quote}

Clearly, the issue of whether, or not, there is a \textit{polynomial time} algorithm $\mathcal{A}$ that gets as input a Boolean formula $f$ and outputs 1 if and only if $f$ is a tautology is meaningful only if we can establish that there \textit{is} an algorithm $\mathcal{A}$ that gets as input a Boolean formula $f$ and outputs 1 if and only if $f$ is a tautology.

Accordingly I show in Section \ref{sec:2.1.1} how it follows from Theorem VII\footnote{cf. \cite{Go31}, p.29: Every recursive relation is arithmetical.} of Kurt G\"{o}del's seminal 1931 paper---on formally undecidable arithmetical propositions---that every recursive function $f(x_{1}, x_{2})$ is representable in PA by a formula $[F(x_{1}, x_{2}, x_{3})]$ such that $[(\exists_{1} x_{3})F(x_{1}, x_{2},$ $x_{3})]$\footnote{The symbol `$[\exists_{1}]$' denotes uniqueness, in the sense that the PA formula $[(\exists_{1} x_{3})F(x_{1}, x_{2}, x_{3})]$ is a short-hand notation for the PA formula $[\neg(\forall x_{3})\neg F(x_{1}, x_{2}, x_{3}) \wedge (\forall y)(\forall z)(F(x_{1}, x_{2}, y) \wedge F(x_{1}, x_{2}, z) \rightarrow y=z)]$.} is ``efficiently verifiable", but not ``efficiently computable", \textit{if} the standard interpretation of PA\footnote{See Section \ref{sec:5.1}. I shall follow Alfred Tarski's terminology and definitions of the satisfaction and truth of the formulas of a formal language under an interpretation as detailed in Section \ref{sec:5}.}---which presumes that Aristotle's particularisation holds over the structure $\mathcal{N}$ of the natural numbers---is sound.

\begin{quote}
\noindent \textbf{Aristotle's particularisation} This holds that from an assertion such as:

\begin{quote}
`It is not the case that, for any given $x$, $P^{*}(x)$\footnote{Notation: The asterisk indicates that the expression is to be interpreted semantically with respect to some well-defined interpretation. I shall aim to use this notation consistently in this investigation.} does not hold',
\end{quote}

\noindent usually denoted symbolically by `$\neg(\forall x)\neg P^{*}(x)$', we may always validly infer in the classical, Aristotlean, logic of predicates\footnote{\cite{HA28}, pp.58-59.} that:

\begin{quote}
`There exists an unspecified $x$ such that $P^{*}(x)$ holds',
\end{quote}

\noindent usually denoted symbolically by `$(\exists x)P^{*}(x)$\footnote{See Appendix A, Section \ref{sec:A} for the meaning and usage of the symbol denoting the existential quantifier in an interpretation.}'.
\end{quote}

However, as L.\ E.\ J.\ Brouwer pointed out in a seminal 1908 paper\footnote{\cite{Br08}.}, the presumption that Aristotle's particularisation holds over $\mathcal{N}$ lies beyond our common intuition. In the rest of the investigation I therefore consider whether the above conclusion would persist under any sound interpretation of PA.

\begin{quote}
\footnotesize{
We may express Aristotle's particularisation in a contemporary context as:

\begin{quote}
\noindent From an assertion such as:

\begin{quote}
`It is not the case that, for any given $x$, any witness\footnote{The word `witness' is intended to be construed broadly in its usual dictionary sense, and not as a specifically defined technical term. However, see Section \ref{sec:5} for a more specific sense of the term `witness' as used in this investigation.} $W_{\mathcal{D}}$ of a domain $\mathcal{D}$ can decide that $P^{*}(x)$ does not hold in $\mathcal{D}$',
\end{quote}

\noindent usually denoted symbolically by `$\neg(\forall x)\neg P^{*}(x)$', we may always validly infer that:

\begin{quote}
`There exists an unspecified $x$ such that any witness $W_{\mathcal{D}}$ of $\mathcal{D}$ can decide that $P^{*}(x)$ holds in $\mathcal{D}$',
\end{quote}

\noindent usually denoted symbolically by `$(\exists x)P^{*}(x)$'.
\end{quote}

I note that, prima facie, Brouwer's objection seems valid since Aristotle's particularisation does not hold if we take $\mathcal{D}$ as the domain of the natural numbers and the witness $W_{\mathcal{N}}$ as a Turing machine, since $P^{*}(x)$ may be a Halting-type of number-theoretic relation. 

Thus, to ensure that the arguments of this investigation are intuitionistically unobjectionable, any assumption that the `standard' interpretation $\mathcal{I}_{PA(\mathcal{N},\ Standard)}$ of PA is sound shall be explicit.
}
\end{quote}

\subsection{Defining instantiational computability and algorithmic computability}
\label{sec:1.02}

We introduce the two concepts\footnote{My thanks to Dr. Chaitanya H. Mehta for advising that the focus of this investigation should be the distinction between these two concepts.}: 

\begin{definition}
\label{sec:1.02.def.1}
\textbf{Instantiational computability}: A Boolean number-theoretic function\footnote{Strictly speaking, a formula of a formal language that interprets as a Boolean number-theoretic function under a well-defined interpretation of the language.} $[F(x_{1}, \ldots, x_{n})]$\footnote{I shall use square brackets to differentiate a formal expression such as $[F(a_{1}, \ldots, a_{n})]$ from its interpretation $F^{*}(x_{1}, \ldots, x_{n})$. See Appendix A, Section \ref{sec:A} for the notation and definitions of standard terms as used in this investigation.} is instantiationally computable if, and only if, there is a Turing machine TM that, for any given sequence of numerals $[(a_{1}, \ldots, a_{n})]$, will accept the natural number input $m$ if $m$ is a unique identification number of the formula $[F(a_{1}, \ldots, a_{n})]$, and will always then halt with one of the following as output:

\begin{quote}
(i) $0$ if $[F(a_{1}, \ldots, a_{n})]$ computes as $0$ (or interprets as true) in $\mathcal{N}$;

\vspace{+.5ex}
(ii) $1$ if $[F(a_{1}, \ldots, a_{n})]$ computes as $1$ (or interprets as false) in $\mathcal{N}$.
\end{quote}
\end{definition}

\begin{definition}
\label{sec:1.02.def.2}
\textbf{Algorithmic computability}: A Boolean number-theoretic function $[F(x_{1},$ $\ldots, x_{n})]$ is algorithmically computable if, and only if, there is a Turing machine TM$_{F}$ that, for any given sequence of numerals $[(a_{1}, \ldots, a_{n})]$, will accept the natural number input $m$ if, and only if, $m$ is a unique identification number of the formula $[F(a_{1}, \ldots, a_{n})]$, and will always then halt with one of the following as output:

\begin{quote}
(i) $0$ if $[F(a_{1}, \ldots, a_{n})]$ computes as $0$ (or interprets as true) in $\mathcal{N}$;

\vspace{+.5ex}
(ii) $1$ if $[F(a_{1}, \ldots, a_{n})]$ computes as $1$ (or interprets as false) in $\mathcal{N}$.
\end{quote}

\begin{quote}
\footnotesize{
The set of identification numbers ${m}$ thus corresponds to the set-theoretically defined language $L$ accepted by TM$_{F}$ in Cook's definition of the class P. 
}
\end{quote}
\end{definition}

It is reasonable to assume that the following thesis will hold when the concepts ``efficiently verifiable" and ``efficiently computable" are formalised in any formal system of Arithmetic:

\begin{thesis}
\label{sec:1.02.lem.1}
(a) A PA formula $[F(x_{1}, \ldots, x_{n})]$ is efficiently verifiable if, and only if, it is instantiationally computable.

\vspace{+.5ex}
(b) If a PA formula $[F(x_{1}, \ldots, x_{n})]$ is efficiently computable then it is algorithmically computable.\hfill $\Box$
\end{thesis}

\vspace{+1ex}
It follows that:

\begin{lemma}
\label{sec:1.02.lem.2}
If Thesis \ref{sec:1.02.lem.1} holds, then P$\neq$NP if a PA formula $[F(x_{1}, \ldots, x_{n})]$ is instantiationally computable, but not algorithmically computable.
\end{lemma}

\vspace{+1ex}
\noindent \textbf{Proof} By Thesis \ref{sec:1.02.lem.1}, if $[F(x_{1}, \ldots, x_{n})]$ is instantiationally computable, then it is efficiently verifiable; whereas if $[F(x_{1}, \ldots, x_{n})]$ is not algorithmically computable, then it is not efficiently computable.\hfill $\Box$

\vspace{+1ex}
Lemma \ref{sec:1.02.lem.2} is intended to highlight the fact that the definition of a tautology \textit{only} requires that a Boolean number-theoretic function $f(x_{1}, \ldots, x_{n})$ be computable instantiationally as always true; unless we presume the Church-Turing Thesis, it does \textit{not} require that $f(x_{1}, \ldots, x_{n})$ be partial recursive, and therefore computable algorithmically as always true.

\begin{quote}
\footnotesize{
I shall argue in Section \ref{sec:7.1} that (as in the case of interpretations of PA in Section \ref{sec:5}) it is an implicit belief in the plausibility of---or informal reliance upon\footnote{See, for instance, \cite{Rg87}, p.21, ``Almost all the proofs in this book will use Church's Thesis to some extent".}---the Church-Turing Thesis that obscures the distinction between `instantiational' computability and `algorithmic' computability.
}
\end{quote}

The question thus arises: Is there a Halting-type PA formula $[F(x_{1}, \ldots, x_{n})]$ that is computable instantiationally, but not algorithmically, as always true under a sound interpretation of PA?

\subsection{Is there a Halting-type tautology?}
\label{sec:1.03}

To place this query in perspective I note that:

\begin{lemma}
\label{sec:1.03.lem.1}
If PA has a sound interpretation $\mathcal{I}_{PA(\mathcal{N},\ Sound)}$ over $\mathcal{N}$, then any PA-provable formula $[F(x_{1}, \ldots, x_{n})]$ is instantiationally computable as always true over $N$ under $\mathcal{I}_{PA(\mathcal{N},\ Sound)}$.
\end{lemma}

\vspace{+1ex}
\noindent \textbf{Proof} G\"{o}del has shown how we can algorithmically assign a unique natural (G\"{o}del) number to each PA formula and to each finite sequence of PA formulas\footnote{\cite{Go31}, p.13.}. G\"{o}del has also shown how we can construct a primitive recursive relation $xBy$\footnote{\cite{Go31}, p.22(45)} that holds if, and only if, $x$ is the G\"{o}del number of a proof sequence in PA, and $y$ is the G\"{o}del number of the last formula of the sequence.

Now, if the PA formula $[F(x_{1}, \ldots, x_{n})]$ is PA-provable then, for any given sequence of numerals $[(a_{1}, \ldots, a_{n})]$, the PA formula $[F(a_{1}, \ldots, a_{n})]$ is PA-provable. Hence $xB \lceil [F(a_{1}, \ldots, a_{n})] \rceil$\footnote{$\lceil [F(a_{1}, \ldots, a_{n})] \rceil$ denotes the G\"{o}del number of $[F(a_{1}, \ldots, a_{n})]$.} always holds for some $x$. Since $xBy$ is recursive, there is a Turing machine TM$_{B}$ that will accept $m$ if $m$ is $\lceil [F(a_{1}, \ldots, a_{n})] \rceil$ and halt with output `provable'.

Since a PA-provable formula is true under $\mathcal{I}_{PA(\mathcal{N},\ Sound)}$, the lemma follows.\hfill $\Box$

\vspace{+1ex}
Although the following argument is informal, a formal proof follows immediately from Section \ref{sec:6.2}, where I show in Theorem \ref{sec:6.2.thm.1} that an algorithmic interpretation $\mathcal{I}_{PA(\mathcal{N},\ Algorithmic)}$ of PA---under which a PA-provable formula is algorithmically computable as always true over $\mathcal{N}$---is sound.

\begin{lemma}
\label{sec:6.03.lem.1}
If PA has a sound interpretation $\mathcal{I}_{PA(\mathcal{N},\ Sound)}$ over $\mathcal{N}$, then any PA-provable formula $[F(x_{1}, \ldots, x_{n})]$ is algorithmically computable as always true over $N$.
\end{lemma}

\vspace{+1ex}
\noindent \textbf{Proof} If a PA formula $[F(x_{1}, \ldots, x_{n})]$ is PA-provable, then there is a finite proof sequence in PA whose last member is $[F(x_{1}, \ldots, x_{n})]$. Under any sound interpretation of PA over $\mathcal{N}$ this sequence must\footnote{For a proof of the necessity see Section \ref{sec:6.3}, Theorem \ref{sec:6.3.thm.1}.} interpret as an algorithm (program) of fixed size that, for any sequence $[(a_{1}, \ldots, a_{n})]$ of PA numerals, decides $[F(a_{1}, \ldots, a_{n})]$ as true. This algorithm defines a Turing machine TM$_{F}$ that, for any natural number sequence $(a_{1}, \ldots, a_{n})$, will:

\begin{quote}
(i) accept the natural number $m$ if, and only if, $m$ is the G\"{o}del number of $[F(a_{1}, \ldots, a_{n})]$;

(ii) halt on any such input $m$ with output `true'. 
 \end{quote} 
 
The lemma follows.\hfill $\Box$

\begin{quote}
\footnotesize{
\begin{lemma}
\label{sec:1.03.lem.1.1}
If PA has a sound interpretation $\mathcal{I}_{PA(\mathcal{N},\ Sound)}$ over $\mathcal{N}$, then we may reasonably assume that any PA-provable formula $[F(x_{1}, \ldots, x_{n})]$ is instantiationally computable as always true over $N$ in polynomial time under $\mathcal{I}_{PA(\mathcal{N},\ Sound)}$.
\end{lemma}

\vspace{+1ex}
\noindent \textbf{Proof} If the PA formula $[F(x_{1}, \ldots, x_{n})]$ is PA-provable then there is a sequence $[S_{1}, S_{2}, \ldots, S_{k}]$ of PA formulas each of which is either a PA-axiom, or a consequence of the PA axioms and the preceding formulas in the sequence by application of some Rule of Deduction of PA, and where $S_{k}$ is $[F(x_{1}, \ldots, x_{n})]$.

Further, for any given sequence of numerals $[a_{1}, \ldots, a_{n}]$, the PA formula $[F(x_{1},$ $\ldots, x_{n}) \rightarrow F(a_{1}, \ldots, a_{n})]$ is PA-provable. Hence there is a proof sequence $[S_{1}, S_{2}, \ldots, S_{k}, S_{k+1}, \ldots, S_{l}]$ such that $S_{l}$ is $[F(a_{1}, \ldots, a_{n})]$.

Now it follows in any system of G\"{o}del numbering such as that defined by G\"{o}del in his seminal 1931 paper on formally undecidable arithmetical propositions\footnote{\cite{Go31}, p.13.} that $\lceil [S_{1}, S_{2}, \ldots, S_{k}, S_{k+1}, \ldots, S_{l}] \rceil = \lceil [S_{1}, S_{2}, \ldots, S_{k}] \rceil \star c(a_{1}, \ldots, a_{n})$, where $c(a_{1}, \ldots, a_{n})$ is a constant whose value is determined by the sequence $(a_{1}, \ldots, a_{n})$.

Now $xB \lceil [F(x_{1}, \ldots, x_{n})] \rceil$ holds for $x = \lceil [S_{1}, S_{2}, \ldots, S_{k}] \rceil$, whilst $xB \lceil [F(a_{1}, \ldots,$ $a_{n})] \rceil$ holds for $x = \lceil [S_{1}, S_{2}, \ldots, S_{k}, S_{k+1}, \ldots, S_{l}] \rceil$. Hence if the Turing machine TM$_{B}$ `computes' $[F(x_{1}, \ldots, x_{n})]$ as `true' in time $t$, then it will `compute' $[F(a_{1}, \ldots, a_{n})]$ as `true' in time $t \star c'(a_{1}, \ldots, a_{n})$ where $c'(a_{1}, \ldots, a_{n})$ is a constant whose value is determined by the sequence $(a_{1}, \ldots, a_{n})$.

It is reasonable to presume that $t$ can be treated as a measure that is representative of the length of the program of the Turing machine TM$_{F}$ defined in Lemma \ref{sec:6.03.lem.1}. The lemma follows.\hfill $\Box$
}
\end{quote}

\vspace{+1ex}
I further note that:

\begin{lemma}
\label{sec:1.03.lem.2}
If PA has a sound interpretation $\mathcal{I}_{PA(\mathcal{N},\ Sound)}$ over $\mathcal{N}$, then there is a PA formula $[F]$ which is instantiationally computable as always true over $N$ under $\mathcal{I}_{PA(\mathcal{N},\ Sound)}$  even though $[F]$ is not PA-provable.
\end{lemma}

\vspace{+1ex}
\noindent \textbf{Proof} G\"{o}del has shown how to construct an arithmetical formula with a single variable---say $[R(x)]$\footnote{G\"{o}del refers to this formula only by its G\"{o}del number $r$ (\cite{Go31}, p.25(12)).}---such that $[R(x)]$ is not PA-provable\footnote{G\"{o}del's aim in \cite{Go31} was to show that $[(\forall x)R(x)]$ is not P-provable; by Generalisation it follows, however, that $[R(x)]$ is also not P-provable.}, but $[R(n)]$ is instantiationally PA-provable for any given PA numeral $[n]$. Hence, for any given numeral $[n]$, the PA formula $xB \lceil [R(n)] \rceil$ must hold for some $x$. The lemma follows.\hfill $\Box$

\vspace{+1ex}
The question arises: Is there a Turing machine TM$_{R}$ that, for any given numeral $[n]$, accepts the natural number input $m$ if, and only if, $m$ is the G\"{o}del number of $[R(n)]$, and halts with output $0$ if $[R(n)]$ is true under $\mathcal{I}_{PA(\mathcal{N},\ Sound)}$, and with output $1$ if $[R(n)]$ is false under $\mathcal{I}_{PA(\mathcal{N},\ Sound)}$?

Obviously there can be no such algorithm if $[R(x)]$ is a Halting-type formula such that there would be some putative G\"{o}del number $\lceil [R(n)] \rceil$ on which any putative Turing machine TM$_{R}$ defined as above cannot output either 0 or 1.

\begin{quote}
\footnotesize{
This could be the case if the definition of the formula in question references---either directly or indirectly---algorithmic computations of some number-theoretic functions over $\mathcal{N}$. Such reference occurs in G\"{o}del's definition of $[R(x)]$, which involves an explicit---and deliberate---self-reference. However it also occurs---albeit implicitly---in G\"{o}del's proof that any recursive Boolean function such as $x_{0}=f(x_{1}, x_{2})$ is representable by a PA formula $[F(x_{0}, x_{1}, x_{2})]$\footnote{\cite{Go31}, p.29, Theorem VII.}. The proof involves defining $[F(x_{0}, x_{1}, x_{2})]$ only by its instantiations. Moreover, for any given numerals $[k,\ m]$, the instantiation $[F(k, m, i)]$ is defined in terms of G\"{o}del's $\beta$-function (see Section \ref{sec:2.2.1})---which is such that $\beta(u_{(k, m)}, v_{(k, m)}, i)$ represents the first $m$ terms, i.e. $f(k, 0), f(k, 1),$ $\ldots , f(k, m)$ of $f(k, x_{2})$. Thus $[F(x_{0}, x_{1}, x_{2})]$ implicitly references the values of a putative $\beta(u_{(x_{1}, x_{2})}, v_{(x_{1}, x_{2})}, i)$ which would represent the sequence $f(x_{1}, 0), f(x_{1},$ $1), \ldots, f(x_{1},$ $x_{2})$ for unspecified $x_{1}, x_{2}$ over $\mathcal{N}$.
}
\end{quote}

The thesis that I shall seek to address formally in this investigation\footnote{See Section \ref{sec:6.3}, Corollary \ref{sec:6.3.cor.1.1}.} is thus:

\begin{thesis}
\label{sec:1.03.ths.1}
Under any sound interpretation of PA, G\"{o}del's $[R(x)]$ is instantiationally computable, but not algorithmically computable, as always true in $\mathcal{N}$.
\end{thesis}

Moreover, I shall seek to show why---as Fortnow appears to suggest---resolving the PvNP problem may not be the major issue; the harder part may be altering our attitudes and beliefs so that we can see what is obstructing such a resolution.

\section{Bridging provability and computability}
\label{sec:2}

In a 1956 letter\footnote{See \cite{Go56} for a translation as provided by Juris Hartmanis.} to John von Neumann, G\"{o}del raised an issue of computational complexity that is commonly accepted as a precursor of the PvNP problem:

\begin{quote}
\footnotesize{
One can obviously easily construct a Turing machine, which for every formula $F$ in first order predicate logic and every natural number $n$, allows one to decide if there is a proof of $F$ of length $n$ (length = number of symbols). Let $\Psi (F,n)$ be the number of steps the machine requires for this and let $\phi(n) = max_{F}\Psi(F,n)$. The question is how fast $\phi(n)$ grows for an optimal machine. One can show that $\phi(n) \geq K.n$. If there really were a machine with $f(n) \approx K.n$ (or even $\approx K.n^{2}$), this would have consequences of the greatest importance. Namely, it would obviously mean that in spite of the undecidability of the Entscheidungsproblem, the mental work of a mathematician concerning Yes-or-No questions could be completely replaced by a machine. After all, one would simply have to choose the natural number $n$ so large that when the machine does not deliver a result, it makes no sense to think more about the problem. Now it seems to me, however, to be completely within the realm of possibility that $\phi(n)$ grows that slowly. Since it seems that $\phi(n) = K.n$ is the only estimation which one can obtain by a generalization of the proof of the undecidability of the Entscheidungsproblem and after all $\phi(n) \approx K.n$ (or $\approx K.n^{2}$) only means that the number of steps as opposed to trial and error can be reduced from $N$ to log $N$ (or (log $N$)$^{2}$). However, such strong reductions appear in other finite problems, for example in the computation of the quadratic residue symbol using repeated application of the law of reciprocity. It would be interesting to know, for instance, the situation concerning the determination of primality of a number and how strongly in general the number of steps in finite combinatorial problems can be reduced with respect to simple exhaustive search.
}
\end{quote}

Clearly issues of computational complexity---such as those raised by G\"{o}del above---are finitary concerns involving number-theoretic functions and relations containing quantification over $\mathcal{N}$ that lie naturally within the domains of:

\begin{quote}
(a) First-order Peano Arithmetic PA, which attempts to capture in a formal language the objective essence of how a human intelligence intuitively \textit{reasons} about number-theoretic predicates, and;

(b) Computability Theory, which attempts to capture in a formal language the objective essence of how a human intelligence intuitively \textit{computes} number-theoretic functions.
\end{quote}

Moreover, since G\"{o}del had already shown in 1931 that every recursive relation can be expressed arithmetically\footnote{\cite{Go31}, Theorem VII, p.31.}, his formulation of the computational complexity of a number-theoretic problem in terms of formal arithmetical provability suggests that we ought to persist in seeking, conversely, an algorithmic interpretation of first-order PA\footnote{Part of the finitary consistency proof for PA sought by Hilbert in his `program' (\cite{Hi30}, pp.485-494).} in Computability Theory, so that any number-theoretic problem can be expressed---and addressed---formally in PA, and its solution, if any, interpreted algorithmically in Computability Theory. I investigate this in detail in Section \ref{sec:5}.

\subsection{G\"{o}del's Theorem V and formally unprovable but interpretively true propositions}
\label{sec:2.1.1}

Now, by G\"{o}del's Theorem V\footnote{\cite{Go31}, p.22.}, every recursive relation $f(x_{1}, \ldots, x_{n})$ can be expressed in PA by a formula $[F(x_{1}, \ldots, x_{n})]$ such that, for any given $n$-tuple of natural numbers $a_{1}, \ldots, a_{n}$:

\begin{quote}
If $f(a_{1}, \ldots, a_{n})$ is true, then PA proves $[F(a_{1}, \ldots, a_{n})]$

If $\neg f(a_{1}, \ldots, a_{n})$ is true, then PA proves $[\neg F(a_{1}, \ldots, a_{n})]$
\end{quote}

G\"{o}del relies only on the above to conclude---in his Theorem VI\footnote{\cite{Go31}, p.24.}---the existence of an arithmetical proposition that is formally unprovable in a Peano Arithmetic, but true under a sound interpretation of the Arithmetic.

However, I now show that it is G\"{o}del's Theorem VII\footnote{\cite{Go31}, p.29.} which---for recursive relations of the form  $x_{0}=\phi(x_{1}, \ldots, x_{n})$ defined by the Recursion Rule\footnote{\cite{Me64}, p.120 \& p.132.}---provides an actual blueprint for the construction of PA formulas that are PA-unprovable, but true under the standard interpretation of PA.

Moreover, I shall show that this G\"{o}delian characteristic is merely a reflection of the fact that, by the instantiational nature of their constructive definition in terms of G\"{o}del's $\beta$-function, such formulas are designed to be instantiationally computable, but not algorithmically computable, under the standard interpretation of PA. 

\subsubsection{Every recursive function is representable in PA}
\label{sec:2.2.1}

I note some standard definitions and results (which implicitly presume\footnote{Such an implicit presumption is seen in G\"{o}del's reference in the statement of his Theorem IX to the negation of a universally quantified formula of the restricted functional calculus as indicative of ``the existence of a counter-example" (\cite{Go31}, p.32).} that the standard interpretation of PA is sound, hence quantifiers are interpreted under the assumption that Aristotle's particularisation is valid over $\mathcal{N}$).

G\"{o}del has defined a primitive recursive function---G\"{o}del's $\beta$-function---as\footnote{cf.\ \cite{Go31}, p.31, Lemma 1; \cite{Me64}, p.131, Proposition 3.21.}:

\begin{quote}
$\beta (x_{1}, x_{2}, x_{3}) = rm(1+(x_{3}+ 1) \star x_{2}, x_{1})$
\end{quote}

\noindent where $rm(x_{1}, x_{2})$ denotes the remainder obtained on dividing $x_{2}$ by $x_{1}$.

G\"{o}del showed that:

\begin{lemma}
\label{sec:2.2.2.lem.1}
For any non-terminating sequence of values $f(x_{1}, 0), f(x_{1}, 1), \ldots$, we can construct natural numbers $b, c$ such that:

\begin{quote}
(i) $j = max(n, f(x_{1}, 0), f(x_{1}, 1), \ldots, f(x_{1}, n))$;

(ii) $c = j$!;

(iii) $\beta(b, c, i) = f(x_{1}, i)$ for $0 \leq i \leq n$.
\end{quote}
\end{lemma}

\vspace{+1ex}
\noindent \textbf{Proof} This is a standard result\footnote{cf.\ \cite{Go31}, p.31, p.31, Lemma 1; \cite{Me64}, p.131, Proposition 3.22.}. We reproduce G\"{o}del's original argument of this critical lemma in an Appendix B, Section \ref{sec:5.4.0}.\hfill $\Box$

\vspace{+1ex}
Now we have the standard definition\footnote{\cite{Me64}, p.118.}:

\begin{definition}
\label{sec:2.2.1.def.1}
A number-theoretic function $f(x_{1}, \ldots, x_{n})$ is said to be representable in PA if, and only if, there is a PA formula $[F(x_{1}, \dots, x_{n+1})]$ with the free variables $[x_{1}, \ldots, x_{n+1}]$, such that, for any given natural numbers $k_{1}, \ldots, k_{n+1}$:

\begin{quote}
(i) if $f(k_{1}, \ldots, k_{n}) = k_{n+1}$ then PA proves: $[F(k_{1}, \ldots, k_{n}, k_{n+1})]$;

(ii) PA proves: $[(\exists_{1} x_{n+1})F(k_{1}, \ldots, k_{n}, x_{n+1})]$.
\end{quote}

The function $f(x_{1}, \ldots, x_{n})$ is said to be strongly representable in PA if we further have that:

\begin{quote}
(iii) PA proves: $[(\exists_{1} x_{n+1})F(x_{1}, \ldots, x_{n}, x_{n+1})]$
\end{quote}
\end{definition}

We then have:

\begin{lemma}
\label{sec:2.2.1.lem.1}
$\beta(x_{1}, x_{2}, x_{3})$ is strongly represented in PA by $[Bt(x_{1}, x_{2}, x_{3}, x_{4})]$, which is defined as follows: 

\begin{quote}
$[(\exists w)(x_{1} = ((1 + (x_{3} + 1)\star x_{2}) \star w + x_{4}) \wedge (x_{4} < 1 + (x_{3} + 1) \star x_{2}))]$.
\end{quote}
\end{lemma}

\vspace{+1ex}
\noindent \textbf{Proof} This is a standard result\footnote{cf. \cite{Me64}, p.131, proposition 3.21.}.\hfill $\Box$

\vspace{+1ex}
G\"{o}del further showed that:

\begin{lemma}
\label{sec:2.2.3}
If $f(x_{1}, x_{2})$ is a recursive function defined by:

\begin{quote}
(i) $f(x_{1}, 0) = g(x_{1})$

(ii) $f(x_{1}, (x_{2}+1)) = h(x_{1}, x_{2}, f(x_{1}, x_{2}))$
\end{quote}

\noindent where $g(x_{1})$ and $h(x_{1}, x_{2}, x_{3})$ are recursive functions of lower rank\footnote{cf.\ \cite{Me64}, p.132; \cite{Go31}, p.30(2).} that are represented in PA by well-formed formulas $[G(x_{1}, x_{2})]$ and $[H(x_{1}, x_{2}, x_{3}, x_{4})]$, then $f(x_{1}, x_{2})$ is represented in PA by the following well-formed formula, denoted by $[F(x_{1}, x_{2}, x_{3})]$:

\begin{quote}
$[(\exists u)(\exists v)(((\exists w)(Bt(u, v, 0, w) \wedge G(x_{1}, w))) \wedge Bt(u, v, x_{2}, x_{3}) \wedge (\forall w)(w< x_{2} \rightarrow (\exists y)(\exists z)(Bt(u, v, w, y) \wedge Bt(u, v, (w+1), z) \wedge H(x_{1}, w, y, z)))]$.
\end{quote}
\end{lemma}

\vspace{+1ex}
\noindent \textbf{Proof} This is a standard result\footnote{cf.\ \cite{Go31}, p.31(2); \cite{Me64}, p.132.}. In view of the significance of this lemma for the resolution of the PvNP problem offered in Lemma \ref{sec:2.3.lem.1} below, we reproduce G\"{o}del's original argument and proof of the lemma in Appendix B, Section \ref{sec:5.4.0}.\hfill $\Box$

\subsection{What does ``$[(\exists_{1} x_{3})F(k, m, x_{3})]$ is provable" assert under the standard interpretation of PA?}

Now, if the PA formula $[F(x_{1}, x_{2}, x_{3})]$ represents in PA the recursive function denoted by $f(x_{1}, x_{2})$ then by definition, for any given numerals $[k, m]$, the formula $[(\exists_{1} x_{3})F(k, m, x_{3})]$ is provable in PA; and true under any sound interpretation of PA. We thus have that:

\begin{lemma} 
\label{sec:2.2.3.lem.1}
If we assume that the standard interpretation of PA is sound, then:

\begin{quote}
``$[(\exists_{1} x_{3})F(k,$ $m, x_{3})]$ is true under the standard interpretation of PA"
\end{quote}

\noindent is the assertion that:

\begin{quote}
Given any natural numbers $k, m$, we can construct natural numbers $t_{(k, m)},$ $u_{(k, m)},$ $v_{(k, m)}$---all functions of $k, m$---such that:

\vspace{+1ex}
\begin{quote}

(a) $\beta(u_{(k, m)}, v_{(k, m)}, 0) = g(k)$;

\vspace{+.5ex}
(b) for all $i<m$, $\beta(u_{(k, m)}, v_{(k, m)}, i) = h(k, i, f(k, i))$;

\vspace{+.5ex}
(c) $\beta(u_{(k, m)}, v_{(k, m)}, m) = t_{(k, m)}$;
 
\end{quote}

\vspace{+1ex}
where $f(x_{1}, x_{2})$, $g(x_{1})$ and $h(x_{1}, x_{2}, x_{3})$ are any recursive functions that are formally represented in PA by $F(x_{1}, x_{2},$ $x_{3}), G(x_{1}, x_{2})$ and $H(x_{1}, x_{2}, x_{3},$ $x_{4})$ respectively such that:

\vspace{+1ex}
\begin{quote}
(i)      $f(k, 0) = g(k)$

\vspace{+.5ex}
(ii)     $f(k, (y+1)) = h(k, y, f(k, y))$ for all $y<m$

\vspace{+.5ex}
(iii)    $g(x_{1})$ and $h(x_{1}, x_{2}, x_{3})$ are recursive functions that are assumed to be of lower rank than $f(x_{1}, x_{2})$.
\end{quote}
\end{quote}
\end{lemma}

\noindent \textbf{Proof} For any given natural numbers $k$ and $m$, if $[F(x_{1}, x_{2}, x_{3})]$ interprets as a well-defined arithmetical relation under the standard interpretation of PA, then a Turing-machine can construct the sequences $f(k, 0), f(k, 1), \ldots ,$ $f(k, m)$ and $\beta(u_{(k, m)}, v_{(k, m)}, 0), \beta(u_{(k, m)}, v_{(k, m)}, 1), \ldots, \beta(u_{(k, m)}, v_{(k, m)}, m)$ and verify the assertion.\hfill $\Box$

\subsubsection{If the standard interpretation of PA is sound, then P$\neq$NP}
\label{sec:2.2.2}

We now see that:

\begin{lemma}
\label{sec:2.3.lem.1}
If the standard interpretation of PA is sound, then $[(\exists_{1} x_{3})F(x_{1}, x_{2},$ $x_{3})]$ is computable instantiationally, but not computable algorithmically, as always true over $\mathcal{N}$.
\end{lemma}

\vspace{+1ex}
\noindent \textbf{Proof} We assume that the standard interpretation of PA is sound (hence we may, for instance, conclude `There is some $x$ such that \ldots' from the assertion `It is not the case that for all $x$ it is not the case that \dots' in the domain $\mathcal {N}$ of the interpretation.). It then follows from Lemma \ref{sec:2.2.3.lem.1} that:

\begin{quote}
(1) $[(\exists_{1} x_{3})F(k, m, x_{3})]$ is PA-provable for any given numerals $[k, m]$. Hence $[(\exists_{1} x_{3})F(k, m, x_{3})]$ is true under the standard interpretation of PA. It then follows from the definition of $[F(x_{1}, x_{2}, x_{3})]$ in Lemma \ref{sec:2.2.3} that, for any given natural numbers $k, m$, we can construct some pair of natural numbers $u_{(k, m)}, v_{(k, m)}$---where $u_{(k, m)}, v_{(k, m)}$ are functions of the given natural numbers $k$ and $m$---such that:

\begin{quote}
(a) $\beta(u_{(k, m)}, v_{(k, m)}, i) = f(k, i)$ for $0 \leq i \leq m$;

\vspace{+.5ex}
(b) $F^{*}(k, m, f(k, m))$ holds in $\mathcal{N}$.
\end{quote}

Since $\beta(x_{1}, x_{2}, x_{3})$ is primitive recursive, $\beta(u_{(k, m)}, v_{(k, m)}, i)$ defines a constructible non-terminating sequence $f'(k, 0), f'(k, 1),$ $\ldots$ for any given natural numbers $k$ and $m$ such that:

\begin{quote}
(c) $f(k, i) = f'(k, i)$ for $0 \leq i \leq m$.
\end{quote}

We can thus define a Turing machine TM$_{\beta(u_{(k, m)}, v_{(k, m)}, i)}$ that will accept the natural number input $g$ if $g$ is the G\"{o}del number of the PA formula $[(\exists_{1} x_{3})F(k, m, x_{3})]$, and then halt with output `true'.

Hence $[(\exists_{1} x_{3})F(x_{1}, x_{2}, x_{3})]$ is computable instantiationally as always true over $\mathcal{N}$ under the standard interpretation of PA.

\vspace{+1ex}
(2) Now, the pair of natural numbers $u_{(x_{1}, x_{2})}, v_{(x_{1}, x_{2})}$ are defined such that:

\begin{quote}
(a) $\beta(u_{(x_{1}, x_{2})}, v_{(x_{1}, x_{2})}, i) = f(x_{1}, i)$ for $0 \leq i \leq x_{2}$;

\vspace{+.5ex}
(b) $F^{*}(x_{1}, x_{2}, f(x_{1}, x_{2}))$ holds in $\mathcal{N}$;
\end{quote}

where $v_{(x_{1}, x_{2})}$ is defined in Lemma \ref{sec:2.2.3} as $j$! (see Lemma \ref{sec:2.2.2}), and:

\begin{quote}
(c) $j = max(n, f(x_{1}, 0), f(x_{1}, 1), \ldots, f(x_{1}, x_{2}))$;

\vspace{+.5ex}
(d) $n$ is the `number' of terms in the sequence $f(x_{1}, 0), f(x_{1}, 1),$ $\ldots, f(x_{1}, x_{2})$.
\end{quote}

\vspace{+1ex}
Since $j$ is not definable for a non-terminating sequence, $\beta(u_{(x_{1}, x_{2})},$ $v_{(x_{1}, x{2})}, i)$ we cannot construct a non-terminating sequence $f'(x_{1}, 0),$ $f'(x_{1}, 1), \ldots$ such that:

\begin{quote}
(e) $f(k, i) = f'(k, i)$ for all $i \geq 0$.
\end{quote}

We cannot thus define a Turing machine TM$_{\beta(u_{(x_{1}, x_{2})}, v_{(x_{1}, x{2})}, i)}$ that will accept the natural number input $g$ if, and only if, $g$ is the G\"{o}del number of the PA formula $[(\exists_{1} x_{3})F(k, m, x_{3})]$, and then halt with output `true'.

Hence $[(\exists_{1} x_{3})F(x_{1}, x_{2}, x_{3})]$ is not computable algorithmically as always true over $\mathcal{N}$ under the standard interpretation of PA.
\end{quote}

\noindent The lemma follows.\hfill $\Box$

\vspace{+1ex}
It follows that:

\begin{theorem}
\label{sec:2.3.thm.1}
If Thesis \ref{sec:1.02.lem.1} holds, and the standard interpretation of PA is sound, then P$\neq$NP.
\end{theorem}

\vspace{+1ex}
\noindent \textbf{Proof} By Lemma \ref{sec:2.3.lem.1}, $[(\exists_{1} x_{3})F(x_{1}, x_{2}, x_{3})]$ is computable instantiationally, but not computable algorithmically, as always true over $\mathcal{N}$. By Lemma \ref{sec:1.02.lem.2}, P$\neq$NP.\hfill $\Box$

\vspace{+1ex}
\begin{quote}
\footnotesize{
A critical issue that I do not address in this investigation is whether the PA formula $[F(x_{1}, x_{2}, x_{3}]$ can be considered to interpret under a sound interpretation of PA as a well-defined predicate since the denumerable sequences $f'(k, 0), f'(k, 1),$ $\ldots , f'(k, m), m_{p}$---where $p>0$, and $m_{p}$ is not equal to $m_{q}$ if $p$ is not equal to $q$---are represented by denumerable, distinctly different, functions $\beta(x_{p_{1}}, x_{p_{2}}, i)$ respectively. There are thus denumerable pairs $(x_{p_{1}}, x_{p_{2}})$ for which $\beta(x_{p_{1}}, x_{p_{2}}, i)$ yields any given sequence $f'(k, 0), f'(k, 1),$ $\ldots , f'(k, m)$.
}
\end{quote}

\section{Is G\"{o}del's undecidable arithmetical proposition a one-off anomaly?}
\label{sec:2.0}

It also follows from the preceding section that:

\begin{corollary}
\label{sec:2.3.cor.2}
If the standard interpretation of PA is sound, then $[(\exists_{1} x_{3})F(x_{1}, x_{2},$ $x_{3})]$ is not PA-provable.
\end{corollary}

\vspace{+1ex}
\noindent \textbf{Proof} If $[(\exists_{1} x_{3})F(x_{1}, x_{2}, x_{3})]$ were PA-provable then, by Lemma \ref{sec:6.03.lem.1}, there would be an algorithm that decides $[(\exists_{1} x_{3})F(x_{1}, x_{2}, x_{3})]$ as always true under the standard interpretation of PA. By Lemma \ref{sec:2.3.lem.1} this is not the case. The corollary follows.\hfill $\Box$

\vspace{+1ex}
Now, the counter-intuitive element in G\"{o}del's conclusions in his 1931 paper has occasionally given rise to the perception that `undecidable' G\"{o}delian propositions are artificially constructed anomalies which are not likely to be encountered in, or have any appreciable significance for, mainstream mathematics. However, this may not be a realistic perception since the first part of G\"{o}del's Theorem VI\footnote{\cite{Go31}, p.25(1).} is merely a special case of the following theorem:

\begin{theorem}
\label{sec:2.0.thm.1}
If PA is consistent and the PA formula $[F(x_{1}, x_{2},$ $x_{3})]$ represents the recursive function $f(x_{1}, x_{2})$, then:

\begin{quote}
(a) $[(\exists_{1} x_{3})F(k, m, x_{3})]$ is PA-provable for any given numerals $[k], [m]$;

(b) $[(\exists_{1} x_{3})F(x_{1}, x_{2}, x_{3})]$ is not PA-provable.
\end{quote}
\end{theorem}

\vspace{+1ex}
\noindent \textbf{Proof} (a) By definition $[(\exists_{1} x_{3})F(k, m, x_{3})]$ is PA-provable for any given numerals $[k], [m]$ if the PA formula $[F(x_{1}, x_{2},$ $x_{3})]$ represents the recursive function $f(x_{1}, x_{2})$.

\vspace{+.5ex}
(b) If $[(\exists_{1} x_{3})F(x_{1}, x_{2}, x_{3})]$ were PA-provable, then it would be algorithmically computable as true under the standard interpretation of PA. By Lemma \ref{sec:2.3.lem.1}, this is not the case. The theorem follows.\hfill $\Box$

\vspace{+1ex}
If, further, we assume that the standard interpretation of PA is sound, then it follows from (a) that the second part of G\"{o}del's Theorem VI\footnote{\cite{Go31}, p.25(2).} is a special case of the following:

\begin{lemma}
\label{sec:2.0.cor.1}
If the standard interpretation of PA is sound, then $[\neg (\exists_{1} x_{3})F(x_{1},$ $x_{2}, x_{3})]$ is not PA-provable.
\end{lemma}

\vspace{+1ex}
\noindent \textbf{Proof} $[\neg (\exists_{1} x_{3})F(x_{1}, x_{2}, x_{3})]$ is an abbreviation of the PA formula:

\vspace{+1ex}
$[(\forall x_{3})\neg F(x_{1}, x_{2}, x_{3}) \wedge (\forall y)(\forall z)(F(x_{1}, x_{2}, y) \wedge F(x_{1}, x_{2}, z) \rightarrow y=z)]$.

\vspace{+1ex}
Under any sound interpretation of PA over $\mathcal{N}$, the latter formula interprets as the arithmetical relation denoted by:

\vspace{+1ex}
$(\forall x_{3})\neg F^{*}(x_{1}, x_{2}, x_{3}) \wedge (\forall y)(\forall z)(F^{*}(x_{1}, x_{2}, y) \wedge F^{*}(x_{1}, x_{2}, z) \rightarrow y=z)$.

\vspace{+1ex}
If the standard interpretation of PA is sound then Aristotle's particularisation holds over $\mathcal{N}$, and this relation can be equivalently denoted by:

\vspace{+1ex}
$\neg (\exists x_{3})F^{*}(x_{1}, x_{2}, x_{3}) \wedge (\forall y)(\forall z)(F^{*}(x_{1}, x_{2}, y) \wedge F^{*}(x_{1}, x_{2}, z) \rightarrow y=z)$.

\vspace{+1ex}
It follows that if $[\neg (\exists_{1} x_{3})F(x_{1}, x_{2}, x_{3})]$ is PA-provable, then $\neg (\exists x_{3})F^{*}(x_{1}, x_{2},$ $x_{3})$ is always true over $\mathcal{N}$.

However, this is false since $(\exists x_{3})F^{*}(x_{1}, x_{2}, x_{3})$ is always true over $\mathcal{N}$ by Definition \ref{sec:2.2.1.def.1}. The lemma follows.\hfill $\Box$ 

\vspace{+1ex}
We thus have:

\begin{corollary}
\label{sec:2.0.cor.2}
If the standard interpretation of PA is sound, then $[(\exists_{1} x_{3})F(k, m, x_{3})]$ is undecidable in PA. \hfill $\Box$
\end{corollary}

\vspace{+1ex}
\noindent \textbf{Proof} The corollary follows from Theorem \ref{sec:2.0.thm.1}(b) and Lemma \ref{sec:2.0.cor.1}.\hfill $\Box$ 

\subsection{The significance of omega-consistency and Hilbert's program}
\label{sec:4}

The significance of Corollary \ref{sec:2.0.cor.2} is that, in order to avoid intuitionistic objections to his reasoning in \cite{Go31}, G\"{o}del did not assume that the standard interpretation of PA is sound. Instead, G\"{o}del introduced the syntactic property of $\omega$-consistency as an explicit assumption in his formal reasoning\footnote{\cite{Go31}, p.23 and p.28.}. G\"{o}del explained at some length\footnote{In his introduction on p.9 of \cite{Go31}.} that his reasons for introducing $\omega$-consistency as an explicit assumption in his formal reasoning was to avoid appealing to the stronger, semantic, concept of classical arithmetical truth---a concept which is implicitly based on an intuitionistically objectionable logic that assumes Aristotle's particularisation is valid over $\mathcal{N}$.

However, I now show that if we assume the standard interpretation of PA is sound, then PA is consistent if, and only if, it is $\omega$-consistent.

\subsubsection{Hilbert's $\omega$-Rule}
\label{sec:4.3}

Assuming that PA has a sound interpretation over $\mathcal{N}$, is it true that:

\begin{quote}
\textbf{Algorithmic $\omega$-Rule}: If it is proved that the PA formula $[F(x)]$ interprets as an arithmetical relation $F^{*}(x)$ that is algorithmically decidable as true for any given natural number $n$, then the PA formula $[(\forall x)F(x)]$ can be admitted as an initial formula (\textit{axiom}) in PA?
\end{quote}

The significance of this query is that, as part of his program for giving mathematical reasoning a finitary foundation, Hilbert\footnote{cf.\ \cite{Hi30}, pp.485-494.} proposed an $\omega$-Rule as a finitary means of extending a Peano Arithmetic to a possible completion (i.e. to logically showing that, given any arithmetical proposition, either the proposition, or its negation, is formally provable from the axioms and rules of inference of the extended Arithmetic):

\begin{quote}
\textbf{Hilbert's $\omega$-Rule}: If it is proved that the PA formula $[F(x)]$ interprets as an arithmetical relation $F^{*}(x)$ that is true for any given natural number $n$, then the PA formula $[(\forall x)F(x)]$ can be admitted as an initial formula (\textit{axiom}) in PA.
\end{quote}

Now, in his 1931 paper---which can, not unreasonably, be seen as the outcome of a presumed attempt to validate Hilbert's $\omega$-rule---G\"{o}del introduced the concept of $\omega$-consistency\footnote{\cite{Go31}, p.23.}, from which it follows that:

\begin{lemma}
If we meta-assume Hilbert's $\omega$-rule for PA, then a consistent PA is necessarily $\omega$-consistent\footnote{However, we cannot similarly conclude from the the Algorithmic $\omega$-Rule that a consistent PA is necessarily $\omega$-consistent.}. \hfill $\Box$
\end{lemma}

\vspace{+1ex}
\noindent \textbf{Proof} If the PA formula $[F(x)]$ interprets as an arithmetical relation $F^{*}(x)$ that is true for any given natural number $n$, and the PA formula $[(\forall x)F(x)]$ can be admitted as an initial formula (\textit{axiom}) in PA, $\neg [(\forall x)F(x)]$ cannot be PA-provable if PA is consistent. The lemma follows.\hfill $\Box$

\vspace{+1ex}
Moreover, it follows from G\"{o}del's 1931 paper that one consequence of assuming Hilbert's $\omega$-Rule is that there must, then, be an undecidable arithmetical proposition\footnote{G\"{o}del constructed an arithmetical proposition $[R(x)]$ and showed that, if a Peano Arithmetic is $\omega$-consistent, then both $[(\forall x)R(x)]$ and $[\neg (\forall x)R(x)]$ are unprovable in the Arithmetic (\cite{Go31}, p.25(1), p.26(2)).}; a further consequence of which is that PA is essentially incomplete.

However, since G\"{o}del's argument in this paper---from which he concludes the existence of an undecidable arithmetical proposition---is based on the weaker (i.e., weaker than assuming Hilbert's $\omega$-rule) premise that a consistent PA can be $\omega$-consistent, the question arises whether an even weaker Algorithmic $\omega$-Rule (which, prima facie, does not imply that a consistent PA is necessarily $\omega$-consistent) can yield a finitary completion for PA as sought by Hilbert, albeit for an $\omega$-\textit{in}consistent PA.

\subsubsection{Aristotle's particularisation and $\omega$-consistency}
\label{sec:4.1}

I shall now argue that these issues are related, and that placing them in an appropriate perspective requires questioning not only the persisting belief that Aristotle's 2000-year old logic of predicates---a critical component of which is Aristotle's particularisation---remains valid even when applied over an infinite domain such as $\mathcal{N}$, but also the basis of Brouwer's denial of the Law of the Excluded Middle following his challenge of the belief in 1908\footnote{\cite{Br08}.}.

Now, we have that:

\begin{lemma}
\label{sec:4.1.lem.1}
If PA is consistent but not $\omega$-consistent, then there is some PA formula $[F(x)]$ such that, under any sound interpretation---say $\mathcal{I}_{PA(\mathcal{N},\ Sound)}$---of PA over $\mathcal{N}$:

\begin{quote}
(i) for any given numeral $[n]$, the PA formula $[F(n)]$ interprets as true under $\mathcal{I}_{PA(\mathcal{N},\ Sound)}$;

(ii) the PA formula $[\neg(\forall x)F(x)]$ interprets as true under $\mathcal{I}_{PA(\mathcal{N},\ Sound)}$.
\end{quote}
\end{lemma}

\vspace{+1ex}
\noindent \textbf{Proof} The lemma follows from the definition of $\omega$-consistency and from Tarski's standard definitions\footnote{\cite{Ta33}; see also \cite{Ho01} for an explanatory exposition. However, for standardisation and convenience of expression, I follow the formal exposition of Tarski's definitions given in \cite{Me64}, p.50.} of the satisfaction, and truth, of the formulas of a formal system such as PA under an interpretation as detailed in Section \ref{sec:5}.\hfill $\Box$

\vspace{+1ex}
Further:

\begin{lemma}
\label{sec:4.1.lem.2}
If the PA formula $[\neg(\forall x)F(x)]$ interprets as true under $\mathcal{I}_{PA (\mathcal{N})}$, then it is not the case that, for any given PA numeral $[n]$, the PA formula $[F(n)]$ interprets as true under $\mathcal{I}_{PA (\mathcal{N})}$.
\end{lemma}

\vspace{+1ex}
\noindent \textbf{Proof} The lemma follows from Tarski's standard definitions of the satisfaction, and truth, of the formulas of a formal system such as PA under an interpretation.\hfill $\Box$

\vspace{+1ex}
It follows that:

\begin{lemma}
\label{sec:4.1.lem.3}
If the interpretation $\mathcal{I}_{PA (\mathcal{N})}$ admits Aristotle's particularisation over $\mathcal{N}$\footnote{As, for instance, in \cite{Me64}, pp.51-52 V(ii).}, and the PA formula $[\neg(\forall x)F(x)]$ interprets as true under $\mathcal{I}_{PA (\mathcal{N})}$, then there is some \textit{unspecified} PA numeral $[m]$ such that the PA formula $[F(m)]$ interprets as false under $\mathcal{I}_{PA (\mathcal{N})}$.
\end{lemma}

\vspace{+1ex}
\noindent \textbf{Proof} The lemma follows from Aristotle's particularisation and Tarski's standard definitions of the satisfaction, and truth, of the formulas of a formal system such as PA under an interpretation.\hfill $\Box$

\vspace{+1ex}
Hence:

\begin{lemma}
\label{sec:4.1.lem.4}
If PA is consistent and Aristotle's particularisation holds over $\mathcal{N}$, there can be no PA formula $[F(x)]$ such that, under any sound interpretation $\mathcal{I}_{PA (\mathcal{N},\ Sound)}$ of PA over $\mathcal{N}$:

\begin{quote}
(i) for any given numeral $[n]$, the PA formula $[F(n)]$ interprets as true under $\mathcal{I}_{PA (\mathcal{N},\ Sound)}$;

(ii) the PA formula $[\neg(\forall x)F(x)]$ interprets as true under $\mathcal{I}_{PA (\mathcal{N},\ Sound)}$.
\end{quote}
\end{lemma}

\vspace{+1ex}
\noindent \textbf{Proof} The lemma follows from the previous two lemma.\hfill $\Box$

\vspace{+1ex}
In other words\footnote{The above argument is made explicit in view of Martin Davis' remark in \cite{Da82}, p.129, that such a proof of $\omega$-consistency may be ``\ldots open to the objection of \textit{circularity}".}:

\begin{corollary}
\label{sec:4.1.lem.5}
If PA is consistent and Aristotle's particularisation holds over $\mathcal{N}$, then PA is $\omega$-consistent.\hfill $\Box$
\end{corollary}

It follows that:

\begin{lemma}
\label{sec:4.1.lem.6}
If Aristotle's particularisation holds over $\mathcal{N}$, then PA is consistent if, and only if, it is $\omega$-consistent.
\end{lemma}

\vspace{+1ex}
\noindent \textbf{Proof} If PA is $\omega$-consistent then, since $[n=n]$ is PA-provable for any given PA numeral $[n]$, we cannot have that $[\neg(\forall x)(x=x)]$ is PA-provable. Since an inconsistent PA proves $[\neg(\forall x)(x=x)]$, an $\omega$-consistent PA cannot be inconsistent.\hfill $\Box$

\begin{quote}
\footnotesize{The arguments of this section and of the preceding Section \ref{sec:2.0} thus show that\footnote{See also \cite{An09b}.} J.\ Barkley Rosser's `extension' of G\"{o}del's argument\footnote{\cite{Ro36}.} succeeds in avoiding an explicit assumption of $\omega$-consistency \textit{only} by implicitly appealing to Aristotle's particularisation.}
\end{quote}

\subsection{Is PA $\omega$-\textit{in}consistent?}
\label{sec:4.2}
Now, it follows from the preceding section that:

\begin{corollary}
\label{sec:4.2.cor.1}
If PA is consistent but not $\omega$-consistent, then Aristotle's particularisation does not hold over $\mathcal{N}$. \hfill $\Box$
\end{corollary}

As the classical, `standard', interpretation of PA---say $\mathcal{I}_{PA(\mathcal{N},\ Standard)}$---appeals to Aristotle's particularisation\footnote{See, for instance, \cite{Me64}, p.107 and p.52(V)(ii).}, it follows that:

\begin{corollary}
\label{sec:4.2.cor.2}
If PA is consistent but not $\omega$-consistent, then the standard interpretation $\mathcal{I}_{PA(\mathcal{N},\ Standard)}$ of PA is not sound, and does not yield a model of PA. \hfill $\Box$
\end{corollary}

Now, formal quantification in computational theory is currently interpreted---as in classical logic\footnote{See \cite{Hi25}, p.382; \cite{HA28}, p.48; \cite{Be59}, pp.178 \& 218.}---so as to admit Aristotle's particularisation over $\mathcal{N}$ as axiomatic\footnote{In the sense of being intuitively obvious. See, for instance, \cite{Da82}, p.xxiv; \cite{Rg87}, p.308 (1)-(4); \cite{EC89}, p.174 (4); \cite{BBJ03}, p.102.}.

However, if Aristotle's particularisation does not hold over $\mathcal{N}$, it would explain to some extent why efforts to resolve the PvNP problem by arguments that appeal to classical Aristotlean logic cannot prevail.

\section{The \textit{implicit} Satisfaction Condition in Tarski's inductive assignment of truth-values under an Interpretation}
\label{sec:5}

I shall now show that a sound algorithmic interpretation---say $\mathcal{I}_{PA(\mathcal{N},\ Algorithmic)}$ ---of PA under which Aristotle's particularisation does not hold over $\mathcal{N}$ is implicit in the `standard' interpretation of PA.

Moreover, the interpretation emerges naturally once we make explicit the precise role of the implicit Satisfaction Condition in Tarski's definitive 1933 paper on the concept of truth in the languages of the deductive sciences\footnote{\cite{Ta33}.}.

Now, Tarski essentially defines\footnote{cf. \cite{Me64}, p.51.}:

\begin{definition}
\label{sec:5.def.1}
If $[A]$ is an atomic formula $[A(x_{1}, x_{2}, \ldots, x_{n})]$ of a formal language S, then a sequence $(a_{1}, a_{2}, \ldots, a_{n})$ in the domain $\mathcal{D}$ of an interpretation $\mathcal{I}_{S(\mathcal{D})}$ of S satisfies $[A]$ if, and only if:

\begin{quote}
(i) $[A(x_{1}, x_{2}, \ldots, x_{n})]$ interprets under $\mathcal{I}_{S(\mathcal{D})}$ as a relation $A^{*}(x_{1}, x_{2}, \ldots, x_{n})$ in $\mathcal{D}$ for a witness $\mathcal{W}_{\mathcal{D}}$ of $\mathcal{D}$;

(ii) $A^{*}(a_{1}, a_{2}, \ldots, a_{n})$ holds in $\mathcal{D}$ if, and only if, $\|$SATCON($\mathcal{I}_{S(\mathcal{D})}$)$\|$ holds for a witness $\mathcal{W}_{\mathcal{D}}$ of $\mathcal{D}$;
\end{quote}

\noindent where $\|$SATCON($\mathcal{I}_{S(\mathcal{D})}$)$\|$ is some Satisfaction Condition that is always decidable by a witness $\mathcal{W}_{\mathcal{D}}$ of $\mathcal{D}$.
\end{definition}

Further, Tarski's analysis shows how we can inductively assign truth values of `satisfaction', `truth', and `falsity' as follows to the compound formulas of  a first-order theory S under the interpretation $\mathcal{I}_{S(\mathcal{D})}$ in terms of \textit{only} the satisfiability of the atomic formulas of S over $\mathcal{D}$\footnote{cf. \cite{Me64}, p.51.}:

\begin{definition}
\label{sec:5.def.2}
A sequence $s$ of $\mathcal{D}$ satisfies $[\neg A]$ under $\mathcal{I}_{S(\mathcal{D})}$ if, and only if, $s$ does not satisfy $[A]$;
\end{definition}

\begin{definition}
\label{sec:5.def.3}
A sequence $s$ of $\mathcal{D}$ satisfies $[A \rightarrow B]$ under $\mathcal{I}_{S(\mathcal{D})}$ if, and only if, either it is not the case that $s$ satisfies $[A]$, or $s$ satisfies $[B]$;
\end{definition}

\begin{definition}
\label{sec:5.def.4}
A sequence $s$ of $\mathcal{D}$ satisfies $[(\forall x_{i})A]$ under $\mathcal{I}_{S(\mathcal{D})}$ if, and only if, given any denumerable sequence $t$ of $\mathcal{D}$ which differs from $s$ in at most the $i$'th component, $t$ satisfies $[A]$;
\end{definition}

\begin{definition}
\label{sec:5.def.5}
A well-formed formula $[A]$ of PA is true under $\mathcal{I}_{S(\mathcal{D})}$ if, and only if, given any denumerable sequence $t$ of $\mathcal{D}$, $t$ satisfies $[A]$;
\end{definition}

\begin{definition}
\label{sec:5.def.6}
A well-formed formula $[A]$ of PA is false under under $\mathcal{I}_{S(\mathcal{D})}$ if, and only if, it is not the case that, given any denumerable sequence $t$ of $\mathcal{D}$, $t$ satisfies $[A]$.
\end{definition}

It then follows that\footnote{cf. \cite{Me64}, pp.51-53.}:

\begin{theorem}
\label{sec:5.thm.1}
(\textit{Satisfaction Theorem}) If, for any interpretation $\mathcal{I}_{S(\mathcal{D})}$ of a first-order theory S, there is a Satisfaction Condition $\|$SATCON($\mathcal{I}_{S(\mathcal{D})}$)$\|$ which holds for a witness $\mathcal{W}_{\mathcal{D}}$ of $\mathcal{D}$, then:

\vspace{+.5ex}
(i) The $\Delta_{0}$ formulas of S are decidable as either true or false over $\mathcal{D}$ under $\mathcal{I}_{S(\mathcal{D})}$;

\vspace{+.5ex}
(ii) If the $\Delta_{n}$ formulas of S are decidable as either true or as false over $\mathcal{D}$ under $\mathcal{I}_{S(\mathcal{D})}$, then so are the $\Delta(n+1)$ formulas of S.
\end{theorem}

\vspace{+1ex}
\noindent \textbf{Proof} It follows from the above definitions that:

\vspace{+1ex}
(a) If, for any given atomic formula $[A(x_{1}, x_{2}, \ldots, x_{n})]$ of S, it is decidable by $\mathcal{W}_{\mathcal{D}}$ whether or not a sequence $(a_{1}, a_{2}, \ldots, a_{n})$ of $\mathcal{D}$ satisfies $[A(x_{1}, x_{2}, \ldots, x_{n})]$ in $\mathcal{D}$ under $\mathcal{I}_{S(\mathcal{D})}$ then, for any given compound formula $[A^{1}(x_{1}, x_{2}, \ldots, x_{n})]$ of S containing any one of the logical constants $\neg, \rightarrow, \forall$, it is decidable by $\mathcal{W}_{\mathcal{D}}$ whether or not the sequence $(a_{1}, a_{2}, \ldots, a_{n})$ of $\mathcal{D}$ satisfies $[A^{1}(x_{1}, x_{2}, \ldots, x_{n})]$ in $\mathcal{D}$ under $\mathcal{I}_{S(\mathcal{D})}$;

\vspace{+1ex}
(b) If, for any given compound formula $[B^{n}(x_{1}, x_{2}, \ldots, x_{n})]$ of S containing $n$ of the logical constants $\neg, \rightarrow, \forall$, it is decidable by $\mathcal{W}_{\mathcal{D}}$ whether or not a sequence $(a_{1}, a_{2}, \ldots, a_{n})$ of $\mathcal{D}$ satisfies $[B^{n}(x_{1}, x_{2}, \ldots, x_{n})]$ in $\mathcal{D}$ under $\mathcal{I}_{S(\mathcal{D})}$ then, for any given compound formula $[B^{(n+1)}(x_{1}, x_{2}, \ldots, x_{n})]$ of S containing $n+1$ of the logical constants $\neg, \rightarrow, \forall$, it is decidable by $\mathcal{W}_{\mathcal{D}}$ whether or not the sequence $(a_{1}, a_{2}, \ldots, a_{n})$ of $\mathcal{D}$ satisfies $[B^{(n+1)}(x_{1}, x_{2}, \ldots, x_{n})]$ in $\mathcal{D}$ under $\mathcal{I}_{S(\mathcal{D})}$;

\vspace{+1ex}
We thus have that:

\vspace{+1ex}
(c) The $\Delta_{0}$ formulas of S are decidable by $\mathcal{W}_{\mathcal{D}}$ as either true or false over $\mathcal{D}$ under $\mathcal{I}_{S(\mathcal{D})}$;

\vspace{+1ex}
(d) If the $\Delta_{n}$ formulas of S are decidable by $\mathcal{W}_{\mathcal{D}}$ as either true or as false over $\mathcal{D}$ under $\mathcal{I}_{S(\mathcal{D})}$, then so are the $\Delta(n+1)$ formulas of S. \hfill $\Box$

\vspace{+1ex}
In other words, if the atomic formulas of of S interpret under $\mathcal{I}_{S(\mathcal{D})}$ as decidable with respect to the Satisfaction Condition $\|$SATCON($\mathcal{I}_{S(\mathcal{D})}$)$\|$ by a witness $\mathcal{W}_{\mathcal{D}}$ over some domain $\mathcal{D}$, then the propositions of S (i.e., the $\Pi_{n}$ and $\Sigma_{n}$ formulas of S) also interpret as decidable with respect to $\|$SATCON($\mathcal{I}_{S(\mathcal{D})}$)$\|$ by the witness $\mathcal{W}_{\mathcal{D}}$ over $\mathcal{D}$.

I now consider the application of Tarski's definitions to various interpretations of first-order Peano Arithmetic PA.

\subsection{The standard interpretation of PA}
\label{sec:5.1}

The standard interpretation $\mathcal{I}_{PA(\mathcal{N},\ Standard)}$ of PA is obtained if, in $\mathcal{I}_{S(\mathcal{D})}$:

\begin{quote}
(a) we define S as PA with standard first-order predicate calculus as the underlying logic\footnote{Where the string $[(\exists \ldots)]$ is defined as---and is to be treated as an abbreviation for---the string $[\neg (\forall \ldots) \neg]$. We do not consider the case where the underlying logic is Hilbert's formalisation of Aristotle's logic of predicates in terms of his $\epsilon$-operator (\cite{Hi27}, pp.465-466).};

(b) we define $\mathcal{D}$ as $\mathcal{N}$;

(c) for any atomic formula $[A(x_{1}, x_{2}, \ldots, x_{n})]$ of PA and sequence $(a_{1}, a_{2}, \ldots, a_{n})$ of $\mathcal{N}$, we take $\|$SATCON($\mathcal{I}_{PA(\mathcal{N})}$)$\|$ as:

\begin{quote}
$\|$$A^{*}(a_{1}^{*}, a_{2}^{*}, \ldots, a_{n}^{*})$ holds in $\mathcal{N}$ and, for any given sequence  $(b_{1}^{*}, b_{2}^{*}, \ldots, b_{n}^{*})$ of $\mathcal{N}$, the proposition $A^{*}(b_{1}^{*}, b_{2}^{*}, \ldots, b_{n}^{*})$ is decidable in $\mathcal{N}$$\|$;
\end{quote}

(d) we define the witness $\mathcal{W}_{(\mathcal{N},\ Standard)}$ informally as the `mathematical intuition' of a human intelligence for whom $\|$SATCON($\mathcal{I}_{PA(\mathcal{N})}$)$\|$ is always \textit{effectively} decidable in $\mathcal{N}$;

\begin{quote}
\footnotesize{
\begin{lemma}
\label{sec:5.1.lem.1}
$A^{*}(x_{1}, x_{2}, \ldots, x_{n})$ is always effectively decidable in $\mathcal{N}$ by $\mathcal{W}_{(\mathcal{N},\ Standard)}$. 
\end{lemma}

\vspace{+1ex}
\textbf{Proof} If $[A(x_{1}, x_{2}, \ldots, x_{n})]$ is an atomic formula of PA then, for any given sequence of numerals $[b_{1}, b_{2}, \ldots, b_{n}]$, the PA formula $[A(b_{1}, b_{2},$ $\ldots, b_{n})]$ is an atomic formula of the form $[c=d]$, where $[c]$ and $[d]$ are atomic PA formulas that denote PA numerals. Since $[c]$ and $[d]$ are recursively defined formulas in the language of PA, it follows from a standard result\footnote{For any natural numbers $m,\ n$, if $m \neq n$, then PA proves $[\neg(m = n)]$ (\cite{Me64}, p.110, Proposition 3.6). The converse is obviously true.} that, if PA is consistent, then $[c=d]$ is algorithmically computable as either true or false in $\mathcal{N}$. In other words, if PA is consistent, then $[A(x_{1}, x_{2}, \ldots, x_{n})]$ is algorithmically computable (since there is a Turing machine TM$_{A}$ that, for any given sequence of numerals $[b_{1}, b_{2}, \ldots, b_{n}]$, will accept the natural number $m$ if, and only if, $m$ is the G\"{o}del number of the PA formula $[A(b_{1}, b_{2},$ $\ldots, b_{n})]$, and halt with output 0 if $[A(b_{1}, b_{2},$ $\ldots, b_{n})]$ interprets as true in $\mathcal{N}$; and halt with output 1 if $[A(b_{1}, b_{2},$ $\ldots, b_{n})]$ interprets as false in $\mathcal{N}$). The lemma follows.\hfill $\Box$
}
\end{quote}

(e) we postulate that Aristotle's particularisation holds over $\mathcal{N}$\footnote{Hence a PA formula such as $[(\exists x)F(x)]$ interprets under $\mathcal{I}_{PA(\mathcal{N},\ Standard)}$ as `There is some natural number $n$ such that $F(n)$ holds in $\mathcal{N}$.}.
\end{quote}

\vspace{+1ex}
Clearly, (e) does not form any part of Tarski's inductive definitions of the satisfaction, and truth, of the formulas of PA under the above interpretation. Moreover, its inclusion makes $\mathcal{I}_{PA(\mathcal{N},\ Standard)}$ extraneously non-finitary\footnote{\cite{Br08}.}.

\begin{quote}
\footnotesize{
The question arises: Can we formulate the `standard' interpretation of PA without assuming (e) extraneously?

I answer this question affirmatively in Section \ref{sec:5.3} where:

\vspace{+.5ex}
(1) I replace the `mathematical intuition' of a human intelligence by defining an `objective' witness $\mathcal{W}_{(\mathcal{N},\ Instantiational)}$ as the meta-theory $\mathcal{M}_{PA}$ of PA;

\vspace{+.5ex}
(2) I show that $\mathcal{W}_{(\mathcal{N},\ Instantiational)}$ can decide whether $c^{*}=d^{*}$ is true or false by instantiationally computing the Boolean function $A^{*}(x_{1}, x_{2}, \ldots, x_{n})$ for any given sequence of natural numbers $(b^{*}_{1}, b^{*}_{2}, \ldots, b^{*}_{n})$;

\vspace{+.5ex}
(3) I show that this yields an instantiational interpretation of PA over $[\mathcal{N}]$ that is sound if, and only if, (e) holds.

\vspace{+.5ex}
(4) $\mathcal{W}_{(\mathcal{N},\ Instantiational)}$ is thus an instantiational formulation of the standard interpretation of PA over $\mathcal{N}$ (which is presumed to be sound).
}
\end{quote}

I note further that if PA is $\omega$-\textit{in}consistent, then Aristotle's particularisation does not hold over $\mathcal{N}$, and the interpretation $\mathcal{I}_{PA(\mathcal{N},\ Standard)}$ is not sound.

\subsection{G\"{o}del's non-standard interpretation of PA}
\label{sec:5.2}

A non-standard (G\"{o}delian) interpretation $\mathcal{I}_{PA(\mathcal{N}_{\omega},\ Non-standard)}$ of a putative $\omega$-consistent PA is obtained if, in $\mathcal{I}_{S(\mathcal{D})}$:

\begin{quote}
(a) we define S as PA with standard first-order predicate calculus as the underlying logic;

(b) we define $\mathcal{D}$ as an undefined extension $\mathcal{N}_{\omega}$ of $\mathcal{N}$;

(c) for any atomic formula $[A(x_{1}, x_{2}, \ldots, x_{n})]$ of PA and sequence $(a_{1}, a_{2}, \ldots, a_{n})$ of $\mathcal{N}_{\omega}$, we take $\|$SATCON($\mathcal{I}_{PA(\mathcal{N}_{\omega})}$)$\|$ as:

\begin{quote}
$\|$$A^{*}(a_{1}^{*}, a_{2}^{*}, \ldots, a_{n}^{*})$ holds in $\mathcal{N}_{\omega}$ and, for any given sequence  $(b_{1}^{*}, b_{2}^{*}, \ldots, b_{n}^{*})$ of $\mathcal{N}_{\omega}$, the proposition $A^{*}(b_{1}^{*}, b_{2}^{*}, \ldots, b_{n}^{*})$ is decidable as either holding or not holding in $\mathcal{N}_{\omega}$$\|$;
\end{quote}

(d) we \textit{postulate} that $\|$SATCON($\mathcal{I}_{PA(\mathcal{N}_{\omega})}$)$\|$ is always decidable by a putative witness $\mathcal{W}_{\mathcal{N_{\omega}}}$, and that $\mathcal{W}_{\mathcal{N_{\omega}}}$ can, further, determine some numbers in $\mathcal{N}_{\omega}$ which are not natural numbers;

(e) we assume that PA is $\omega$-consistent.
\end{quote}

Clearly, the interpretation $\mathcal{I}_{PA(\mathcal{N}_{\omega},\ Non-standard)}$ of a putative $\omega$-consistent PA cannot claim to be finitary. Moreover, if PA is $\omega$-\textit{in}consistent, then the G\"{o}delian non-standard interpretation $\mathcal{I}_{PA(\mathcal{N}_{\omega},\ Non-standard)}$ of PA is also not sound\footnote{In which case we cannot validly conclude from G\"{o}del's formal reasoning in (\cite{Go31}) that PA must have a non-standard model.}.

\subsection{An instantiational interpretation of PA in PA}
\label{sec:5.3}

I next consider the instantiational interpretation $\mathcal{I}_{PA(\mathcal{N},\ Instantiational)}$ of PA where:

\begin{quote}
(a) we define S as PA with standard first-order predicate calculus as the underlying logic;

(b) we define $\mathcal{D}$ as PA;

(c) for any atomic formula $[A(x_{1}, x_{2}, \ldots, x_{n})]$ of PA and any sequence $[(a_{1}, a_{2}, \ldots, a_{n})]$ of PA numerals, we take $\|$SATCON($\mathcal{I}_{PA(\mathcal{PA})}$)$\|$ as:

\begin{quote}
$\|$$[A(a_{1}, a_{2}, \ldots, a_{n})]$ is provable in PA and, for any given sequence of numerals $[(b_{1}, b_{2}, \ldots, b_{n})]$ of PA, the formula $[A(b_{1}, b_{2}, \ldots, b_{n})]$ is decidable as either provable or not provable in PA$\|$;
\end{quote}

(d) we define the witness $\mathcal{W}_{(\mathcal{N},\ Instantiational)}$ as the meta-theory $\mathcal{M}_{PA}$ of PA.

\begin{quote}
\footnotesize{
\begin{lemma}
\label{sec:5.3.lem.1}
$[A(x_{1}, x_{2}, \ldots, x_{n})]$ is always effectively decidable in PA by $\mathcal{W}_{(\mathcal{N},\ Instantiational)}$. 
\end{lemma}

\vspace{+1ex}
\textbf{Proof} It follows from G\"{o}del's definition of the primitive recursive relation $xBy$\footnote{\cite{Go31}, p. 22(45).}---where $x$ is the G\"{o}del number of a proof sequence in PA whose last term is the PA formula with G\"{o}del-number $y$---that, if $[A(x_{1}, x_{2}, \ldots, x_{n})]$ is an atomic formula of PA, $\mathcal{M}_{PA}$ can effectively decide instantiationally for any given sequence $[(b_{1}, b_{2}, \ldots, b_{n})]$ of PA numerals which one of the PA formulas $[A(b_{1}, b_{2}, \ldots, b_{n})]$ and $[\neg A(b_{1}, b_{2}, \ldots, b_{n})]$ is necessarily PA-provable.\hfill $\Box$
}
\end{quote}
\end{quote}

Now, if PA is consistent but not $\omega$-consistent, then there is a G\"{o}delian formula $[R(x)]$ such that (see Section \ref{sec:6.3}):

\begin{quote}
(i) $[(\forall x)R(x)]$ is not PA-provable;

(ii) $[\neg (\forall x)R(x)]$ is PA-provable;

(iii) for any given numeral $[n]$, $[R(n)]$ is PA-provable.
\end{quote}

However, if $\mathcal{I}_{PA(\mathcal{N},\ Instantiational)}$ is sound, then (ii) implies contradictorily that it is not the case that, for any given numeral $[n]$, $[R(n)]$ is PA-provable.

\vspace{+1ex}
It follows that if $\mathcal{I}_{PA(\mathcal{N},\ Instantiational)}$ is sound then PA is $\omega$-consistent and, ipso facto, Aristotle's particularisation must hold over $\mathcal{N}$.

Moreover, if PA is consistent, then every PA-provable formula interprets as true under some sound interpretation of PA. Hence $\mathcal{M}_{PA}$ can effectively decide whether, for any given sequence of natural numbers $(b_{1}^{*}, b_{2}^{*},$ $\ldots, b_{n}^{*})$ in $\mathcal{N}$, the proposition $A^{*}(b_{1}^{*}, b_{2}^{*}, \ldots, b_{n}^{*})$ holds or not in $\mathcal{N}$.

It follows that $\mathcal{I}_{PA(\mathcal{N},\ Instantiational)}$ is an instantiational formulation of the `standard' interpretation of PA in which we do not need to extraneously assume that Aristotle's particularisation holds over $\mathcal{N}$.

The interpretation $\mathcal{I}_{PA(\mathcal{N},\ Instantiational)}$ is of interest because, if it were a sound interpretation of PA, then PA would establish its own consistency\footnote{cf. G\"{o}del's Theorem XI in \cite{Go31}, p.36.}!

\subsection{A set-theoretic interpretation of PA}
\label{sec:5.4}

I consider next a set-theoretic interpretation $\mathcal{I}_{PA(\mathcal{ZF},\ Cantor)}$ of PA over the domain of ZF sets, which is obtained if, in $\mathcal{I}_{S(\mathcal{D})}$:

\begin{quote}
(a) we define S as PA with standard first-order predicate calculus as the underlying logic;

(b) we define $\mathcal{D}$ as ZF;

(c) for any atomic formula $[A(x_{1}, x_{2}, \ldots, x_{n})]$ and sequence $[(a_{1}, a_{2}, \ldots, a_{n})]$ of PA, we take $\|$SATCON($\mathcal{I}_{S(\mathcal{ZF})}$)$\|$ as:

\begin{quote}
$\|$$[A^{*}(a_{1}^{*}, a_{2}^{*}, \ldots, a_{n}^{*})]$ is provable in ZF and, for any given sequence  $[(b_{1}^{*}, b_{2}^{*}, \ldots, b_{n}^{*})]$ of ZF, the formula $[A^{*}(b_{1}^{*}, b_{2}^{*}, \ldots, b_{n}^{*})]$ is decidable as either provable or not provable in ZF$\|$;
\end{quote}

(d) we define the witness $\mathcal{W}_{\mathcal{ZF}}$ as the meta-theory $\mathcal{M}_{ZF}$ of ZF which can always decide effectively whether or not $\|$SATCON($\mathcal{I}_{S(\mathcal{ZF})}$)$\|$ holds in ZF.
\end{quote}

Now, if the set-theoretic interpretation $\mathcal{I}_{PA(\mathcal{ZF},\ Cantor)}$ of PA is sound, then every sound interpretation of ZF would, ipso facto, be a sound interpretation of PA. In Appendix C, Section \ref{sec:5.4.1} I show, however, that this is not the case, and so the set-theoretic interpretation $\mathcal{I}_{PA(\mathcal{ZF},\ Cantor)}$ of PA is not sound.

\subsection{A purely algorithmic interpretation of PA}
\label{sec:5.5}

I finally consider the purely algorithmic interpretation $\mathcal{I}_{PA(\mathcal{N},\ Algorithmic)}$ of PA where:

\begin{quote}
(a) we define S as PA with standard first-order predicate calculus as the underlying logic;

(b) we define $\mathcal{D}$ as $\mathcal{N}$;

(c) for any atomic formula $[A(x_{1}, x_{2}, \ldots, x_{n})]$ and sequence $(a_{1}, a_{2}, \ldots, a_{n})$ of $\mathcal{N}$, we take $\|$SATCON($\mathcal{I}_{PA(\mathcal{N})}$)$\|$ as:

\begin{quote}
$\|$$A^{*}(a_{1}^{*}, a_{2}^{*}, \ldots, a_{n}^{*})$ holds in $\mathcal{N}$ and, for any given sequence  $(b_{1}^{*}, b_{2}^{*}, \ldots, b_{n}^{*})$ of $\mathcal{N}$, the proposition $A^{*}(b_{1}^{*}, b_{2}^{*}, \ldots, b_{n}^{*})$ is decidable as either holding or not holding in $\mathcal{N}$$\|$;
\end{quote}

(d) we define the witness $\mathcal{W}_{(\mathcal{N},\ Algorithmic)}$ as a Turing machine TM$_{A^{*}}$ for whom $\|$SATCON($\mathcal{I}_{PA(\mathcal{N})}$)$\|$ is always \textit{effectively} decidable in $\mathcal{N}$:
\end{quote}

\begin{quote}
\footnotesize{
\begin{lemma}
\label{sec:5.5.lem.1}
$A^{*}(x_{1}, x_{2}, \ldots, x_{n})$ is always effectively decidable in $\mathcal{N}$ by $\mathcal{W}_{(\mathcal{N},\ Al-}$ $_{gorithmic)}$. 
\end{lemma}

\vspace{+1ex}
\textbf{Proof} If $[A(x_{1}, x_{2}, \ldots, x_{n})]$ is an atomic formula of PA then, for any given sequence of numerals $[b_{1}, b_{2}, \ldots, b_{n}]$, the PA formula $[A(b_{1}, b_{2},$ $\ldots, b_{n})]$ is an atomic formula of the form $[c=d]$, where $[c]$ and $[d]$ are atomic PA formulas that denote PA numerals. Since $[c]$ and $[d]$ are recursively defined formulas in the language of PA, it follows from a standard result\footnote{For any natural numbers $m,\ n$, if $m \neq n$, then PA proves $[\neg(m = n)]$ (\cite{Me64}, p.110, Proposition 3.6). The converse is obviously true.} that, if PA is consistent, then $[c=d]$ is algorithmically computable as either true or false in $\mathcal{N}$. In other words, if PA is consistent, then $[A(x_{1}, x_{2}, \ldots, x_{n})]$ is algorithmically computable (since there is a Turing machine TM$_{A}$ that, for any given sequence of numerals $[b_{1}, b_{2}, \ldots, b_{n}]$, will accept the natural number $m$ if, and only if, $m$ is the G\"{o}del number of the PA formula $[A(b_{1}, b_{2},$ $\ldots, b_{n})]$, and halt with output 0 if $[A(b_{1}, b_{2},$ $\ldots, b_{n})]$ interprets as true in $\mathcal{N}$; and halt with output 1 if $[A(b_{1}, b_{2},$ $\ldots, b_{n})]$ interprets as false in $\mathcal{N}$). The lemma follows.\hfill $\Box$
}
\end{quote}

It follows that $\mathcal{I}_{PA(\mathcal{N},\ Algorithmic)}$ is an algorithmic formulation of the `standard' interpretation of PA in which we do not extraneously assume that Aristotle's particularisation holds over $\mathcal{N}$.

\begin{quote}
\footnotesize{
I shall show that if $\mathcal{I}_{PA(\mathcal{N},\ Algorithmic)}$ is sound, then PA is \textit{not} $\omega$-consistent. Hence Aristotle's particularisation does not hold over $\mathcal{N}$, and the interpretation is finitary and intuitionistically unobjectionable. Moreover---since the Law of the Excluded Middle is provable in an $\omega$-inconsistent PA (and therefore holds in $\mathcal{N}$)---it achieves this without the discomforting, stringent, Intuitionistic requirement that we reject the underlying logic of PA!
}
\end{quote}

I now show that $\mathcal{I}_{PA(\mathcal{N},\ Algorithmic)}$ is sound, and consider the consequences for the PvNP problem and for Church's Thesis.

\section{The algorithmic interpretation of PA is sound}
\label{sec:6}

In Section \ref{sec:5} of this investigation I defined the two interpretations $\mathcal{I}_{PA(\mathcal{N},\ Standard)}$ and $\mathcal{I}_{PA(\mathcal{N},\ Algorithmic)}$ in terms of Tarski's\footnote{\cite{Ta33}.} inductive definitions of the satisfaction, and truth, of the formulas of a formal system under an interpretation. It thus follows by induction on $k$ that\footnote{cf. \cite{Me64}, pp.51-53.}:

\subsection{Interpreting quantification}
\label{sec:6.1}

\begin{lemma}
\label{sec:6.1.lem.1}
(\textit{Universal:Standard}) A $\Pi_{k}$ PA formula such as $[(\forall x)A(x)]$ interprets as true\footnote{See Definition \ref{sec:5.def.5}} under $\mathcal{I}_{PA(\mathcal{N},\ Standard)}$ if, and only if, for any given natural number $n$, $A^{*}(n)$ is true in $\mathcal{N}$.
\end{lemma}

\vspace{+1ex}
\noindent \textbf{Proof} The lemma follows from the definition of $\mathcal{I}_{PA(\mathcal{N},\ Standard)}$ by induction on $k$.\hfill $\Box$

\begin{lemma}
\label{sec:6.1.lem.2}
(\textit{Universal:Algorithmic}) A $\Pi_{k}$ PA formula such as $[(\forall x)A(x)]$ interprets as true under $\mathcal{I}_{PA(\mathcal{N},\ Algorithmic)}$ if, and only if,
$A^{*}(x)$ is algorithmically computable as always true in $\mathcal{N}$.
\end{lemma}

\vspace{+1ex}
\noindent \textbf{Proof} The lemma follows from the definition of $\mathcal{I}_{PA(\mathcal{N},\ Algorithmic)}$ by induction on $k$.\hfill $\Box$

\begin{lemma}
\label{sec:6.1.lem.3}
(\textit{Existential:Standard}) A $\Sigma_{k}$ PA formula such as $[(\exists x)A(x)]$\footnote{Note that $[(\exists x)A(x)]$ is merely the abbreviation for $[\neg(\forall x)\neg A(x)]$.} interprets as true under $\mathcal{I}_{PA(\mathcal{N},\ Standard)}$ if, and only if, it is not true that, for any given natural number $n$, $A^{*}(n)$ is false in $\mathcal{N}$, \textbf{and} we may conclude that there exists some natural number $n$ such that $A^{*}(n)$ holds in $\mathcal{N}$\footnote{Since Aristotle's particularisation is assumed to hold in $\mathcal{N}$ under $\mathcal{I}_{PA(\mathcal{N},\ Standard)}$.}.
\end{lemma}

\vspace{+1ex}
\noindent \textbf{Proof} The lemma follows from the definition of $\mathcal{I}_{PA(\mathcal{N},\ Standard)}$ by induction on $k$.\hfill $\Box$

\begin{lemma}
\label{sec:6.1.lem.4}
(\textit{Existential:Algorithmic}) A $\Sigma_{k}$ PA formula such as $[(\exists x)$ $A(x)]$ interprets as true under $\mathcal{I}_{PA(\mathcal{N},\ Algorithmic)}$ if, and only if, $\neg A^{*}(x)$ is not algorithmically computable as always true in $\mathcal{N}$, \textbf{but} we may \textbf{not} conclude that there exists some natural number $n$ such that $A^{*}(n)$ holds in $\mathcal{N}$\footnote{Since $A^{*}(x)$ may be a Halting-type of relation such that, for any given natural number $n$, it is meta-mathematically---even if not algorithmically---decidable that $A^{*}(n)$ is false. As I show in Section \ref{sec:6.3}, G\"{o}del's relation $R(x)$ is precisely such a relation.}.
\end{lemma}

\vspace{+1ex}
\noindent \textbf{Proof} The lemma follows from the definition of $\mathcal{I}_{PA(\mathcal{N},\ Algorithmic)}$ by induction on $k$.\hfill $\Box$

\subsection{Interpreting the PA axioms}
\label{sec:6.2}

We note first that:

\begin{lemma}
\label{sec:6.2.lem.1}
The PA axioms PA$_{1}$ to PA$_{8}$ are algorithmically computable as always true over $\mathcal{N}$ under the interpretation $\mathcal{I}_{PA(\mathcal{N},\ Algorithmic)}$.
\end{lemma}

\vspace{+1ex}
\noindent \textbf{Proof} Since $[x+y]$, $[x \star y]$, $[x = y]$, $[{x^{\prime}}]$ are defined recursively\footnote{cf. \cite{Go31}, p.17.}, the PA axioms PA$_{1}$ to PA$_{8}$ interpret as recursive relations that do not involve any quantification. The lemma follows.\hfill $\Box$

\vspace{+1ex}
\noindent Further:

\begin{lemma}
\label{sec:6.2.lem.2}
For any given PA formula $[F(x)]$, the Induction axiom schema $[F(0) \rightarrow (((\forall x)(F(x) \rightarrow F(x^{\prime}))) \rightarrow (\forall x)F(x))]$ interprets as true under $\mathcal{I}_{PA(\mathcal{N},}$ $_{\ Algorithmic)}$.
\end{lemma}

\vspace{+1ex}
\noindent \textbf{Proof} By Tarski's Definitions \ref{sec:5.def.1} to \ref{sec:5.def.6}:

\begin{quote}

(a) If $[F(0)]$ interprets as false under $\mathcal{I}_{PA(\mathcal{N},}$ $_{\ Algorithmic)}$ the lemma is proved.

(b) If $[F(0)]$ interprets as true and $[(\forall x)(F(x) \rightarrow F(x^{\prime}))]$ interprets as false under $\mathcal{I}_{PA(\mathcal{N},}$ $_{\ Algorithmic)}$, the lemma is proved.

(c) If $[F(0)]$ and $[(\forall x)(F(x) \rightarrow F(x^{\prime}))]$ both interpret as true under $\mathcal{I}_{PA(\mathcal{N},}$ $_{\ Algorithmic)}$, then by the Satisfaction Theorem \ref{sec:5.thm.1} and the algorithmic interpretation $\mathcal{I}_{PA(\mathcal{N},}$ $_{\ Algorithmic)}$ defined in Section \ref{sec:5.5}, $[F(x) \rightarrow F(x^{\prime})]$ is algorithmically computable as always true over $\mathcal{N}$ under $\mathcal{I}_{PA(\mathcal{N},}$ $_{\ Algorithmic)}$.

There is thus a Turing machine TM$_{F}$ such that, for any natural number $n$, TM$_{F}$ will accept the natural number $m$ if, and only if, $m$ is the G\"{o}del number of $[F(n) \rightarrow F(n^{\prime})]$ and will halt with output `true'.

Since $[F(0)]$ interprets as true under $\mathcal{I}_{PA(\mathcal{N},}$ $_{\ Algorithmic)}$, it follows that there is a Turing machine TM$_{F^{\prime}}$ such that, for any natural number $n$, TM$_{F^{\prime}}$ will accept the natural number $m$ if, and only if, $m$ is the G\"{o}del number of $[F(n)]$ and will halt with output `true'.

Hence $[(\forall x)F(x)]$ is algorithmically computable as always true under $\mathcal{I}_{PA(\mathcal{N},}$ $_{\ Algorithmic)}$.
\end{quote}

Since the above cases are exhaustive, the lemma follows.\hfill $\Box$

\begin{quote}
\footnotesize{
I note that the interpretation $\mathcal{I}_{PA(\mathcal{N},\ Algorithmic)}$ settles the Poincar\'{e}-Hilbert debate\footnote{See \cite{Hi27}, p.472; also \cite{Br13}, p.59; \cite{We27}, p482; \cite{Pa71}, p.502-503.} in the latter's favour. Poincar\'{e} believed that the Induction Axiom could not be justified finitarily, as any such argument would necessarily need to appeal to infinite induction. Hilbert believed that a finitary proof of the consistency of PA was possible.
}
\end{quote}

\vspace{+1ex}
\noindent Similarly:

\begin{lemma}
\label{sec:6.2.lem.3}
Generalisation preserves truth under $\mathcal{I}_{PA(\mathcal{N},\ Algorithmic)}$.
\end{lemma}

\vspace{+1ex}
\noindent \textbf{Proof} The two meta-assertions:

\begin{quote}
`$[F(x)]$ interprets as true under $\mathcal{I}_{PA(\mathcal{N},}$ $_{\ Algorithmic)}$\footnote{See Definition \ref{sec:5.def.5}}'

and

`$[(\forall x)F(x)]$ interprets as true under $\mathcal{I}_{PA(\mathcal{N},\ Algorithmic)}$'
\end{quote}

\noindent both mean:

\begin{quote}
$[F(x)]$ is algorithmically computable as always true under $\mathcal{I}_{PA(\mathcal{N},}$ $_{\ Algorithmic)}$.\hfill $\Box$
\end{quote} 

\vspace{+1ex}
It is also straightforward to see that:

\begin{lemma}
\label{sec:6.2.lem.4}
Modus Ponens preserves truth under $\mathcal{I}_{PA(\mathcal{N},\ Algorithmic)}$. \hfill $\Box$
\end{lemma}

We thus have that:

\begin{lemma}
\label{sec:6.2.lem.5}
The axioms of PA are always true under the finitary interpretation $\mathcal{I}_{PA(\mathcal{N},\ Algorithmic)}$, and the rules of inference of PA preserve the properties of satisfaction/truth under $\mathcal{I}_{PA(\mathcal{N},\ Algorithmic)}$.\hfill $\Box$
\end{lemma}

Hence:

\begin{theorem}
\label{sec:6.2.thm.1}
The interpretation $\mathcal{I}_{PA(\mathcal{N},\ Algorithmic)}$ of PA is sound.
\end{theorem}

We thus have a finitary proof that:

\begin{theorem}
\label{sec:6.2.thm.2}
PA is consistent. \hfill $\Box$
\end{theorem}

\section{A Provability Theorem for PA}
\label{sec:6.3}

I now show that PA can have no non-standard model, since it is `algorithmically' complete in the sense that:

\begin{theorem}
\label{sec:6.3.thm.1}
(Provability Theorem for PA) A PA formula $[F(x)]$ is PA-provable if, and only if, $[F(x)]$ is algorithmically computable in $\mathcal{N}$.
\end{theorem}

\vspace{+1ex}
\noindent \textbf{Proof} We have by definition that $[(\forall x)F(x)]$ interprets as true under the interpretation $\mathcal{I}_{PA(\mathcal{N},\ Algorithmic)}$ if, and only if, $[F(x)]$ is algorithmically computable in $\mathcal{N}$.

Since $\mathcal{I}_{PA(\mathcal{N},\ Algorithmic)}$ is sound, it defines a finitary model of PA over $\mathcal{N}$---say $\mathcal{M}_{PA(\beta)}$---such that:

\begin{quote}
If $[(\forall x)F(x)]$ is PA-provable, then $[F(x)]$ is algorithmically computable in $\mathcal{N}$;

If $[\neg(\forall x)F(x)]$ is PA-provable, then it is not the case that $[F(x)]$ is algorithmically computable in $\mathcal{N}$.
\end{quote}

Now, we cannot have that both $[(\forall x)F(x)]$ and $[\neg(\forall x)F(x)]$ are PA-unprovable for some PA formula $[F(x)]$, as this would yield the contradiction:

\begin{quote}
(i) There is a finitary model---say $M1_{\beta}$---of PA+$[(\forall x)F(x)]$ in which $[F(x)]$ is algorithmically computable in $\mathcal{N}$.

(ii) There is a finitary model---say $M2_{\beta}$---of PA+$[\neg(\forall x)F(x)]$ in which it is not the case that $[F(x)]$ is algorithmically computable in $\mathcal{N}$.
\end{quote}

The lemma follows.\hfill $\Box$

\begin{corollary}
PA is categorical.
\end{corollary}

\vspace{+1ex}
By the argument in Theorem \ref{sec:6.3.thm.1} it follows that:

\begin{corollary}
\label{sec:6.3.cor.1}
The PA formula $[\neg(\forall x)R(x)]$ defined in Lemma \ref{sec:1.03.lem.2} is PA-provable. \hfill $\Box$
\end{corollary}

\begin{corollary}
\label{sec:6.3.cor.1.1}
Under any sound interpretation of PA, G\"{o}del's $[R(x)]$ interprets as an instantiationally computable, but not algorithmically computable, tautology over $N$.
\end{corollary}

\vspace{+1ex}
\noindent \textbf{Proof} G\"{o}del has shown that $[R(x)]$\footnote{G\"{o}del refers to this formula only by its G\"{o}del number $r$; \cite{Go31}, p.25, eqn.12.} interprets as an instantiationally computable tautology\footnote{\cite{Go31}, p.26(2): ``$(n)\neg(nB_{\kappa}(17Gen\ r))$ holds"}. By Corollary \ref{sec:6.3.cor.1} $[R(x)]$ is not algorithmically computable as always true in $\mathcal{N}$.\hfill $\Box$

\begin{theorem}
\label{sec:6.3.thm.2}
P$\neq$NP.
\end{theorem}

\vspace{+1ex}
\noindent \textbf{Proof} By Corollary \ref{sec:6.3.cor.1.1}, $[R(x)]$ is instantiationally computable, but not algorithmically computable over $\mathcal{N}$. The theorem follows immediately from Lemma \ref{sec:1.02.lem.2}.\hfill $\Box$

\begin{corollary}
\label{sec:6.3.cor.2}
PA is \textit{not} $\omega$-consistent.\footnote{This conclusion is contrary to accepted dogma. See, for instance, Davis' remarks in \cite{Da82}, p.129(iii) that ``\ldots there is no equivocation. Either an adequate arithmetical logic is $\omega$-inconsistent (in which case it is possible to prove false statements within it) or it has an unsolvable decision problem and is subject to the limitations of G\"{o}del's incompleteness theorem".}
\end{corollary}

\vspace{+1ex}
\noindent \textbf{Proof} G\"{o}del has shown that if PA is consistent, then $[R(n)]$ is PA-provable for any given PA numeral $[n]$\footnote{\cite{Go31}, p.26(2).}. By Corollary \ref{sec:6.3.cor.1} and the definition of $\omega$-consistency, if PA is consistent then it is \textit{not} $\omega$-consistent.\hfill $\Box$

\begin{corollary}
\label{sec:6.3.cor.3}
The standard interpretation $\mathcal{I}_{PA(\mathcal{N},\ Standard)}$ of PA is not sound, and does not yield a model of PA\footnote{I note that finitists of all hues---ranging from Brouwer \cite{Br08} to Alexander Yessenin-Volpin \cite{He04}---have persistently questioned the soundness of the `standard' interpretation $\mathcal{I}_{PA(\mathcal{N},\ Standard)}$.}.
\end{corollary}

\vspace{+1ex}
\noindent \textbf{Proof} By Corollary \ref{sec:4.2.cor.1} if PA is consistent but not $\omega$-consistent, then Aristotle's particularisation does not hold over $\mathcal{N}$. Since the `standard', interpretation of PA appeals to Aristotle's particularisation, the lemma follows.\hfill $\Box$

\begin{quote}
\footnotesize{
Since formal quantification is currently interpreted in classical logic\footnote{See \cite{Hi25}, p.382; \cite{HA28}, p.48; \cite{Be59}, pp.178 \& 218.} so as to admit Aristotle's particularisation over $\mathcal{N}$ as axiomatic\footnote{In the sense of being intuitively obvious. See, for instance, \cite{Da82}, p.xxiv; \cite{Rg87}, p.308 (1)-(4); \cite{EC89}, p.174 (4); \cite{BBJ03}, p.102.}, the above suggests that we may need to review number-theoretic arguments\footnote{For instance---as shown in Sections \ref{sec:2.0} and \ref{sec:4}---Rosser's construction of an undecidable arithmetical proposition in PA (see \cite{Ro36})---which does not explicitly assume that PA is $\omega$-consistent---implicitly presumes that Aristotle's particularisation holds over $\mathcal{N}$.} that appeal unrestrictedly to classical Aristotlean logic.
}
\end{quote}

\subsection{The Provability Theorem for PA and Bounded Arithmetic}
\label{sec:6.3.1}

In a 1997 paper\footnote{\cite{Bu97}.}, Samuel R. Buss considered Bounded Arithmetics obtained by:

\begin{quote}
(a) limiting the applicability of the Induction Axiom Schema in PA only to functions with quantifiers bounded by an unspecified natural number bound $b$;

(b) `weakening' the statement of the axiom with the aim of differentiating between effective computability over the sequence of natural numbers, and feasible `polynomial-time' computability over a bounded sequence of the natural numbers\footnote{See also \cite{Pa71}.}.
\end{quote}

Presumably Buss' intent---as expressed below---is to build a bridge between provability in a Bounded Arithmetic and Computability so that a $\Pi_{k}$ formula, say $[(\forall x)f(x)]$, is provable in the Bounded Arithmetic if, and only if, there is an algorithm that, for any given numeral $[n]$, decides the $\Delta_{(k/(k-1))}$ formula $[f(n)]$ as `true':

\begin{quote}
If $[(\forall x)(\exists y)f(x, y)]$ is provable, then there should be an algorithm to find $y$ as a function of $x$\footnote{See \cite{Bu97}.}.
\end{quote}

Since we have proven such a Provability Theorem for PA in the previous section, the first question arises:

\vspace{+1ex}
\noindent \textbf{Does the introduction of bounded quantifiers yield any computational advantage?}

\vspace{+1ex}
Now, one difference\footnote{I suspect the only one.} between a Bounded Arithmetic and PA is that we can presume in the Bounded Arithmetic that, from a proof of $[(\exists y)f(n, y)]$, we may always conclude that there is some numeral $[m]$ such that $[f(n, m)]$ is provable in the arithmetic; however, this is not a sound conclusion in PA.

Reason: Since $[(\exists y)f(n, y)]$ is simply a shorthand for $[\neg (\forall y)\neg f(n, y)]$, such a presumption implies that Aristotle's particularisation holds over the natural numbers under any sound interpretation of PA.

To see that (as Brouwer steadfastly held) this may not always be the case, interpret $[(\forall x)f(x)]$ as\footnote{We have seen in the earlier sections that such an interpretation is sound.}:

\begin{quote}
There is an algorithm that decides $[f(n)]$ as `true' for any given numeral $[n]$.
\end{quote}

In such case, if $[(\forall x)(\exists y)f(x, y)]$ is provable in PA, then we can only conclude that:

\begin{quote}
There is an algorithm that, for any given numeral $[n]$, decides that it is not the case that there is an algorithm that, for any given numeral $[m]$, decides $[\neg f(n, m)]$ as `true'.
\end{quote}

We cannot, however, conclude - as we can in a Bounded Arithmetic - that:

\begin{quote}
There is an algorithm that, for any given numeral $[n]$, decides that there is an algorithm that, for some numeral $[m]$, decides $[f(n, m)]$ as `true'.
\end{quote}

Reason: $[(\exists y)f(n, y)]$ may be a Halting-type formula for some numeral $[n]$.

This could be the case if $[(\forall x)(\exists y)f(x, y)]$ were PA-\textit{un}provable, but $[(\exists y)f(n, y)]$ PA-provable for any given numeral $[n]$.

Presumably it is the belief that any sound interpretation of PA requires Aristotle's particularisation to hold in $\mathcal{N}$, and the recognition that the latter does not admit linking provability to computability in PA, which has led to considering the effect of bounding quantification in PA.

However, as we have seen in the preceding sections, we are able to link provability to computability through the Provability Theorem for PA by recognising precisely that, to the contrary, any interpretation of PA which requires Aristotle's particularisation to hold in $\mathcal{N}$ cannot be sound!

The postulation of an unspecified bound in a Bounded Arithmetic in order to arrive at a provability-computability link thus appears dispensible.

The question then arises:

\vspace{+1ex}
\noindent \textbf{Does `weakening' the PA Induction Axiom Schema yield any computational advantage?}

\vspace{+1ex}
Now, Buss considers a bounded arithmetic $S_{2}$ which is, essentially, PA with the following `weakened' Induction Axiom Schema, PIND\footnote{Where $\lfloor \frac{x}{2} \rfloor$ denotes the largest natural number lower bound of the rational $\frac{x}{2}$.}:

\vspace{+1ex}
\noindent $[\{f(0)\ \&\ (\forall x)(f(\lfloor \frac{x}{2} \rfloor) \rightarrow f(x))\} \rightarrow (\forall x)f(x)]$

\vspace{+1ex}
However, PIND can be expressed in first-order Peano Arithmetic PA as follows:

\vspace{+1ex}
\noindent $[\{f(0)\ \&\ (\forall x)(f(x) \rightarrow (f(2*x)\ \&\ f(2*x+1)))\} \rightarrow (\forall x)f(x)]$.

\vspace{+1ex}
Moreover, the above is a particular case of PIND($k$):

\vspace{+1ex}
\noindent $[\{f(0)\ \&\ (\forall x)(f(x) \rightarrow (f(k*x)\ \&\ f(k*x+1)\ \&\ \ldots \& f(k*x+k-1)))\} \rightarrow (\forall x)f(x)]$.

\vspace{+1ex}
Now we have the PA theorem: 

\vspace{+1ex}
\noindent $[(\forall x)f(x) \rightarrow \{f(0)\ \&\ (\forall x)(f(x) \rightarrow f(x+1))\}]$

\vspace{+1ex}
It follows that the following is also a PA theorem:

\vspace{+1ex}
\noindent $[\{f(0)\ \&\ (\forall x)(f(x) \rightarrow f(x+1))\} \rightarrow \{f(0)\ \&\ (\forall x)(f(x) \rightarrow (f(k*x)\ \&\ f(k*x+1)\ \&\ \ldots \&\ f(k*x+k-1)))\}]$

\vspace{+1ex}
In other words, for any numeral $[k]$, PIND($k$) is equivalent in PA to the standard Induction Axiom of PA!

Thus, the Provability Theorem for PA suggests that all arguments and conclusions of a Bounded Arithmetic can be reflected in PA without any loss of generality.

\section{Church's Thesis is false}
\label{sec:7.1}

One reason why efforts to prove P$=$NP remain unsuccessful may lie in recognising that---contrary to accepted dogma\footnote{See \cite{Kl52}, p.300; \cite{Me64}, p.227; \cite{Rg87}, p.20; \cite{EC89}, p.85; \cite{BBJ03}, p23.}---the term `effective computability' \textit{can} be precisely defined. For instance, we can define:

\begin{definition}
\label{sec:7.1.def.3}
A number-theoretic function is effectively computable if, and only if, it is computable instantiationally.
\end{definition}

Prima facie, this definition adequately captures our intuitive understanding of the term `effective computability'.

Now, classical theory argues that (standard results):

\begin{lemma}
\label{sec:7.1.lem.1}
Every Turing-computable function (or relation, treated as a Boolean function) $F$ is partial recursive , and, if $F$ is total , then $F$ is recursive\footnote{cf.\ \cite{Me64}, p.233, Corollary 5.13.}. \hfill $\Box$
\end{lemma}

\begin{lemma}
\label{sec:7.1.lem.2}
Every partial recursive function (or relation, treated as a Boolean function) is Turing-computable\footnote{cf.\ \cite{Me64}, p.237, Corollary 5.15.}. \hfill $\Box$
\end{lemma}

It follows that the following---essentially unverifiable but refutable---theses are classically equivalent\footnote{cf.\ \cite{Me64}, p.237.}:

\begin{quote}
\textbf{Standard Church's Thesis}\footnote{\textit{Church's (original) Thesis} The effectively computable number-theoretic functions are the algorithmically computable number-theoretic functions \cite{Ch36}.} A number-theoretic function (or relation, treated as a Boolean function) is effectively computable if, and only if, it is partial-recursive\footnote{cf.\ \cite{Me64}, p.227.}.

\textbf{Standard Turing's Thesis}\footnote{After describing what he meant by ``computable" numbers in the opening sentence of his 1936 paper on Computable Numbers \cite{Tu36}, Turing immediately expressed this thesis---albeit informally---as: ``\ldots the computable numbers include all numbers which could naturally be regarded as computable".} A number-theoretic function (or relation, treated as a Boolean function) is effectively computable if, and only if, it is Turing-computable\footnote{cf.\ \cite{BBJ03}, p.33.}.
\end{quote}

\subsection{The Church and Turing Theses are false}
\label{sec:7.2}

However, Church's Thesis ignores the doctrine of Occam's razor by postulating a strong identity---and not simply a weak equivalence---between an effectively computable number-theoretic function and some algorithmically computable function.

Consequently, Church's Thesis (Turing's Thesis) does not admit the possibility of an arithmetical function $F$ that is computable instantiationally but not algorithmically. It follows that:

\begin{theorem}
\label{sec:7.2.thm.1}
The Church and Turing theses do not hold.
\end{theorem}

\vspace{+1ex}
\noindent \textbf{Proof} By Corollary \ref{sec:6.3.cor.1.1} G\"{o}del's $[R(x)]$ is instantiationally computable as always true, but it is not algorithmically computable as always true. The lemma follows.\hfill $\Box$

\subsection{Recognising instantiational computability as `effective'}
\label{sec:7.3}

It is significant that G\"{o}del (initially) and Church (subsequently---possibly under the influence of G\"{o}del's disquietitude) enunciated Church's formulation of `effective computability' as a Thesis because G\"{o}del was instinctively uncomfortable with accepting it as a definition that fully captures the essence of `\textit{intuitive} effective computability'\footnote{See \cite{Si97}.}.

G\"{o}del's reservations seem vindicated if we accept that a number-theoretic function can be computable instantiationally, but not algorithmically.

The possibility that `truth' may be be `effectively' decidable instantiationally, but not algorithmically, is implicit in G\"{o}del's famous 1951 Gibbs lecture\footnote{\cite{Go51}.}, where he remarks: 

\begin{quote}
``I wish to point out that one may conjecture the truth of a universal proposition (for example, that I shall be able to verify a certain property for any integer given to me) and at the same time conjecture that no general proof for this fact exists. It is easy to imagine situations in which both these conjectures would be very well founded. For the first half of it, this would, for example, be the case if the proposition in question were some equation $F(n) = G(n)$ of two number-theoretical functions which could be verified up to very great numbers $n$."\footnote{Parikh's paper \cite{Pa71} can also be viewed as an attempt to investigate the consequences of expressing the essence of G\"{o}del's remarks formally.} 
\end{quote}

Such a possibility is also implicit in Turing's remarks\footnote{\cite{Tu36}, \S9, para II.}:

\begin{quote}
``The computable numbers do not include all (in the ordinary sense) definable numbers. Let P be a sequence whose \textit{n}-th figure is 1 or 0 according as \textit{n} is or is not satisfactory. It is an immediate consequence of the theorem of \S8 that P is not computable. It is (so far as we know at present) possible that any assigned number of figures of P can be calculated, but not by a uniform process. When sufficiently many figures of P have been calculated, an essentially new method is necessary in order to obtain more figures." 
\end{quote}

The need for placing such a distinction on a formal basis has also been expressed explicitly on occasion\footnote{Parikh's distinction between `decidability' and `feasibility' in \cite{Pa71} also appears to echo the need for such a distinction.}. Thus, Boolos, Burgess and Jeffrey\footnote{\cite{BBJ03}, p. 37.} define a diagonal function, $d$, any value of which can be decided effectively, although there is no single algorithm that can effectively compute $d$. 

Now, the straightforward way of expressing this phenomenon should be to say that there are well-defined number-theoretic functions that are instantiationally computable, but not algorithmically computable\footnote{Or, preferably, one could borrow the analogous terminology from the theory of functions of real and complex variables and term such functions as computable, but not uniformly computable.}. Yet, following Church and Turing, such functions are labeled as effectively uncomputable\footnote{The issue here seems to be that, when using language to express the abstract objects of our individual, and common, mental `concept spaces', we use the word `exists' loosely in three senses, without making explicit distinctions between them (see \cite{An07}).}!

\begin{quote}
\noindent ``According to Turing's Thesis, since $d$ is not Turing-computable, $d$ cannot be effectively computable. Why not? After all, although no Turing machine computes the function $d$, we were able to compute at least its first few values, For since, as we have noted, $f_{1} = f_{1} = f_{1} =$ the empty function we have $d(1) = d(2) = d(3) = 1$. And it may seem that we can actually compute $d(n)$ for any positive integer $n$---if we don't run out of time."\footnote{\cite{BBJ03}, p.37.}
\end{quote}

The reluctance to treat a function such as $d(n)$---or the function $\Omega(n)$ that computes the $n^{th}$ digit in the decimal expression of a Chaitin constant $\Omega$\footnote{Chaitin's Halting Probability is given by $0 < \Omega = \sum2^{-|p|} < 1$, where the summation is over all self-delimiting programs $p$ that halt, and $|p|$ is the size in bits of the halting program $p$; see \cite{Ct75}.}---as computable, on the grounds that the `time' needed to compute it increases monotonically with $n$, is curious\footnote{The incongruity of this is addressed by Parikh in \cite{Pa71}.}; the same applies to any total Turing-computable function $f(n)$\footnote{The only difference being that, in the latter case, we know there is a common `program' of constant length that will compute $f(n)$ for any given natural number $n$; in the former, we know we may need distinctly different programs for computing $f(n)$ for different values of $n$, where the length of the program will, sometime, reference $n$.}!

\section{Conclusions}
\label{sec:1.04}

I have defined what it means for a number-theoretic function to be:

\begin{quote}
(i) Instantiationally computable;

\vspace{+.5ex}
(ii) Algorithmically computable.
\end{quote}

I then show that:

\begin{quote}
\begin{conclusion} If Aristotle's particularisation is presumed valid over the structure $\mathcal{N}$ of the natural numbers---as is the case under the standard interpretation of PA---then it follows from the instantiational nature of the constructive definition of the G\"{o}del $\beta$-function that a primitive recursive relation can be instantiationally equivalent to an arithmetical relation where the former is algorithmically computable as always true over $\mathcal{N}$ whilst the latter is instantiationally computable, but not algorithmically computable, as always true over $\mathcal{N}$.
\end{conclusion}
\end{quote}

I then show that:

\begin{quote}
\begin{conclusion}
$\mathcal{I}_{PA(\mathcal{N},\ Standard)}$ is sound if, and only if, Aristotle's particularisation holds over $\mathcal{N}$; and the latter is the case if, and only if, PA is $\omega$-consistent.
\end{conclusion}
\end{quote}

Under the standard interpretation $\mathcal{I}_{PA(\mathcal{N},\ Standard)}$ of PA over the domain $\mathcal{N}$, if $[A]$ is an atomic formula $[A(x_{1}, x_{2}, \ldots, x_{n})]$ of PA, then the sequence of natural numbers $(a_{1}, a_{2}, \ldots, a_{n})$ satisfies $[A]$ if, and only if $A^{*}(a_{1}, a_{2}, \ldots, a_{n})$ holds in $\mathcal{N}$ \textit{and} we presume that Aristotle's particularisation is valid over $\mathcal{N}$.

\vspace{+1ex}
I have now shown that:

\begin{quote}
\begin{conclusion}
We can define a sound interpretation $\mathcal{I}_{PA(\mathcal{N},\ Algorithmic)}$ of PA over the domain $\mathcal{N}$ where, if $[A]$ is an atomic formula $[A(x_{1}, x_{2}, \ldots,$ $x_{n})]$ of PA, then the sequence of natural numbers $(a_{1}, a_{2}, \ldots, a_{n})$ satisfies $[A]$ if, and only if $[A(a_{1}, a_{2}, \ldots, a_{n})]$ is algorithmically computable under $\mathcal{I}_{PA(\mathcal{N},\ Algorithmic)}$, but we do not presume that Aristotle's particularisation is valid over $\mathcal{N}$.
\end{conclusion}
\end{quote}

It follows that:

\begin{quote}
\begin{conclusion}
PA is consistent.
\end{conclusion}
\end{quote}

I have then shown that:

\begin{quote}
\begin{conclusion}
(Provability Theorem for PA) A PA formula $[F(x)]$ is PA-provable if, and only if, $[F(x)]$ is algorithmically computable in $\mathcal{N}$.
\end{conclusion}
\end{quote}

It follows that:

\begin{quote}
\begin{conclusion}
PA is categorical.
\end{conclusion}
\end{quote}

In his 1931 paper, G\"{o}del showed how to construct a formula $[R(x)]$ with a single free variable in any Peano Arithmetic such that:

\begin{quote}
If the first-order Peano Arithmetic PA is assumed to be consistent, then:

\begin{quote}
(i) for any PA numeral $[n]$, the PA formula $[R(n)]$ is provable in PA;

(ii) the PA formula $[R(x)]$ is not provable in PA.
\end{quote}
\end{quote}

I have now shown that:

\begin{quote}
\begin{conclusion}
\label{que:1}
(a) We may conclude from G\"{o}del's argument that $[R(x)]$ is not algorithmically computable as true in $\mathcal{N}$.

(b) $[\neg(\forall x)R(x)]$ is PA-provable;

(c) PA is \textit{not} $\omega$-consistent;

(d) The `standard' interpretation $\mathcal{I}_{PA(\mathcal{N},\ Standard)}$ of PA over $\mathcal{N}$ is not sound.
\end{conclusion}
\end{quote}

Since the above imples that G\"{o}del's formula $[R(x)]$ is computable instantiationally but not algorithmically as true over $\mathcal{N}$, I conclude that:

\begin{quote}
\begin{conclusion}
(i) P$\neq$NP;

(ii) the Church and Turing Theses do not hold?
\end{conclusion}
\end{quote}

\section{Appendix A: Notation, Definitions and Comments}
\label{sec:A}

\noindent \textbf{Notation} In this investigation I use square brackets to indicate that the contents represent a symbol or a formula of a formal theory, generally assumed to be well-formed unless otherwise indicated by the context.

\begin{quote}
\footnotesize{
In other words, expressions inside the square brackets are to be only viewed syntactically as juxtaposition of symbols that are to be formed and manipulated upon strictly in accordance with specific rules for such formation and manipulation---in the manner of a mechanical or electronic device---without any regards to what the symbolism might represent semantically under an interpretation that gives them meaning. 

Moreover, even though the formula `$[F(x)]$' of a formal Arithmetic may interpret as the arithmetical relation expressed by `$F^{*}(x)$', the formula `$[(\exists x)R(x)]$' need not interpret as the arithmetical proposition denoted by the usual abbreviation `$(\exists x)R^{*}(x)$.' The latter denotes the phrase `There is some $x$ such that $R^{*}(x)$'. As Brouwer had noted\footnote{\cite{Br08}; see also \cite{An08}.}, this concept is not always capable of an unambiguous meaning that can be represented in a formal language by the formula `$[(\exists x)R(x)]$' which, in a formal language, is merely an abbreviation for the formula `$[\neg(\forall x)\neg R(x)]$'.

By `expressed' I mean here that the symbolism is simply a short-hand abbreviation for referring to abstract concepts that may, or may not, be capable of a precise `meaning'. Amongst these are symbolic abbreviations which are intended to express the abstract concepts---particularly those of `existence'---involved in propositions that refer to non-terminating processes and infinite aggregates.
}
\end{quote}

\noindent \textbf{Provability} A formula $[F]$ of a formal system S is provable in S (S-provable) if, and only if, there is a finite sequence of S-formulas $[F_{1}], [F_{2}], \ldots, [F_{n}]$ such that $[F_{n}]$ is $[F]$ and, for all $1 \leq i \leq n$, $[F_{i}]$ is either an axiom of S or a consequence of the axioms of S, and the formulas preceding it in the sequence, by means of the rules of deduction of S.

\vspace{+1ex}
\noindent \textbf{The structure $\mathcal{N}$} The structure of the natural numbers---namely, \{$N$ (\textit{the set of natural numbers}); $=$ (\textit{equality}); $'$ (\textit{the successor function}); $+$ (\textit{the addition function}); $ \ast $ (\textit{the product function}); $0$ (\textit{the null} \textit{element})\}.

\vspace{+1ex}
\noindent \textbf{The axioms of first-order Peano Arithmetic (PA)}

\begin{tabbing}
\textbf{PA$_{1}$}  \= $[(x_{1} = x_{2}) \rightarrow ((x_{1} = x_{3}) \rightarrow (x_{2} = x_{3}))]$; \\

\textbf{PA$_{2}$}  \> $[(x_{1} = x_{2}) \rightarrow (x_{1}^{\prime} = x_{2}^{\prime})]$; \\

\textbf{PA$_{3}$}  \> $[0 \neq x_{1}^{\prime}]$; \\

\textbf{PA$_{4}$}  \> $[(x_{1}^{\prime} = x_{2}^{\prime}) \rightarrow (x_{1} = x_{2})]$; \\

\textbf{PA$_{5}$}  \> $[( x_{1} + 0) = x_{1}]$; \\

\textbf{PA$_{6}$}  \> $[(x_{1} + x_{2}^{\prime}) = (x_{1} + x_{2})^{\prime}]$; \\

\textbf{PA$_{7}$}  \> $[( x_{1} \star 0) = 0]$; \\

\textbf{PA$_{8}$}  \> $[( x_{1} \star x_{2}^{\prime}) = ((x_{1} \star x_{2}) + x_{1})]$; \\

\textbf{PA$_{9}$}  \> For any well-formed formula $[F(x)]$ of PA: \\

\> $[F(0) \rightarrow (((\forall x)(F(x) \rightarrow F(x^{\prime}))) \rightarrow (\forall x)F(x))]$.
\end{tabbing}

\vspace{+1ex}
\noindent \textbf{Generalisation in PA} If $[A]$ is PA-provable, then so is $[(\forall x)A]$.

\vspace{+1ex}
\noindent \textbf{Modus Ponens in PA} If $[A]$ and $[A \rightarrow B]$ are PA-provable, then so is $[B]$.

\vspace{+1ex}
\noindent \textbf{Standard interpretation of PA} The standard interpretation $\mathcal{I}_{PA(\mathcal{N},\ Standard)}$ of PA over the structure $\mathcal{N}$ is the one in which the logical constants have their `usual' interpretations\footnote{See \cite{Me64}, p.49.} in Aristotle's logic of predicates\footnote{Thus, Aristotle's particularisation holds over $\mathcal{N}$ in the standard interpretation of PA.}, and\footnote{See \cite{Me64}, p.107.}:

\begin{tabbing}
(a) \= the set of non-negative integers is the domain; \\

(b) \> the integer 0 is the interpretation of the symbol [0]; \\

(c) \> the successor operation (addition of 1) is the interpretation of the $[']$ \\ \> function; \\

(d) \> ordinary addition and multiplication are the interpretations of $[+]$ and \\ \> $[*]$; \\

(e) \> the interpretation of the predicate letter $[=]$ is the identity relation.
\end{tabbing}

\noindent \textbf{Simple consistency} A formal system S is simply consistent if, and only if, there is no S-formula $[F(x)]$ for which both $[(\forall x)F(x)]$ and $[\neg(\forall x)F(x)]$ are S-provable.

\vspace{+1ex}
\noindent \textbf{$\omega$-consistency} A formal system S is $\omega$-consistent if, and only if, there is no S-formula $[F(x)]$ for which, first, $[\neg(\forall x)F(x)]$ is S-provable and, second, $[F(a)]$ is S-provable for any given S-term $[a]$.

\vspace{+.5ex}
\noindent \textbf{Soundness (formal system)}  A formal system S is sound under an interpretation $\mathcal{I}_{S}$ if, and only if, every theorem $[T]$ of S translates as `$[T]$ is true under $\mathcal{I}_{S}$'.

\vspace{+.5ex}
\noindent \textbf{Soundness (interpretation)}  An interpretation $\mathcal{\mathcal{I}_{S}}$ of a formal system S is sound if, and only if, S is sound under the interpretation $\mathcal{\mathcal{I}_{S}}$.

\begin{quote}
\scriptsize{\textbf{Soundness in classical logic}. In classical logic, a formal system $S$ is sometimes defined as `sound' if, and only if, it has an interpretation; and an interpretation is defined as the assignment of meanings to the symbols, and truth-values to the sentences, of the formal system. Moreover, any such interpretation is a model of the formal system. This definition suffers, however, from an implicit circularity: the formal logic $L$ underlying any interpretation of $S$ is implicitly assumed to be `sound'. The above definitions seek to avoid this implicit circularity by delinking the defined `soundness' of a formal system under an interpretation from the implicit `soundness' of the formal logic underlying the interpretation. This admits the case where, even if $L_{1}$ and $L_{2}$ are implicitly assumed to be sound, $S+L_{1}$ is sound, but $S+L_{2}$ is not. Moreover, an interpretation of $S$ is now a model for $S$ if, and only if, it is sound.\footnote{My thanks to Professor Rohit Parikh for highlighting the need for making such a distinction explicit.}}
\end{quote}

\vspace{+1ex}
\noindent \textbf{Categoricity} A formal system S is categorical if, and only if, it has a sound interpretation and any two sound interpretations of S are isomorphic.\footnote{Compare \cite{Me64}, p.91.}

\section{Appendix B: G\"{o}del's Theorem VII(2)}
\label{sec:5.4.0}
\begin{center}
\footnotesize{(Excerpted from \cite{Go31} pp.29-31.)}
\end{center}

\vspace{+1ex}
Every relation of the form $x_{0} = \phi(x_{1}, \ldots, x_{n})$, where $\phi$ is recursive, is arithmetical and we apply complete induction on the rank of $\phi$. Let $\phi$ have rank $s (s>1)$. \ldots

\vspace{+1ex}
$\phi(0, x_{2}, \ldots, x_{n}) = \psi(x_{2}, \ldots, x_{n})$

$\phi(k+1, x_{2}, \ldots, x_{n}) = \mu[k, \phi(k, x_{2}, \ldots, x_{n}), x_{2}, \ldots, x_{n}]$

\noindent(where $\psi, \mu$ have lower rank than $s$).

\vspace{+1ex}
\ldots we apply the following procedure: one can express the relation $x_{0} = \phi(x_{1}, \ldots, x_{n})$ with the help of the concept ``sequence of numbers" $(f)$\footnote{$f$ denotes here a variable whose domain is the sequence of natural numbers. The $(k+1)$st term of a sequence $f$ is designated $f_{k}$ (and the first, $f_{0}$).} in the following manner:

\begin{quote}
$x_{0} = \phi(x_{1}, \ldots, x_{n}) \sim (\exists f)\{f_{0} = \psi(x_{2}, \ldots, x_{n})\ \&\ (\forall k)(k<x_{1} \rightarrow f_{k+1} = \mu(k, f_{k}, x_{2}, \ldots, x_{n})\ \&\ x_{0} = f_{x_{1}}\}$
\end{quote}

If $S(y, x_{2}, \ldots, x_{n}), T(z, x_{1}, \ldots, x_{n+1})$ are the arithmetical relations which, according to the inductive hypothesis, are equivalent to $y = \psi(x_{2}, \ldots, x_{n})$, and $z = \mu(x_{1}, \ldots, x_{n+1})$ respectively, then we have:

\begin{quote}
$x_{0} = \phi(x_{1}, \ldots, x_{n}) \sim (\exists f)\{S(f_{0}, x_{2}, \ldots, x_{n})\ \&\ (\forall k)[k<x_{1} \rightarrow T(f_{k+1}, k, x_{2},$ $\ldots, x_{n})]\ \&\ x_{0} = f_{x_{1}}\}$ \hspace{+35ex} (17)
\end{quote}

Now we replace the concept ``sequence of numbers" by ``pairs of numbers" by correlating with the number pair $n, d$ the sequence of numbers $f^{(n, d)}$ ($f_{k}^{(n, d)} = [n]_{1+(k+1)d}$, where $[n]_{p}$ denotes the smallest non-negative remainder of $n$ modulo $p$).

Then:

\vspace{+1ex}
Lemma 1: If $f$ is an arbitrary sequence of natural numbers and $k$ is an arbitrary natural number, then there exists a pair of natural numbers $n, d$ such that $f^{(n, d)}$ and $f$ coincide in their first $k$ terms.

\vspace{+1ex}
Proof: Let $l$ be the greatest of the numbers $k, f_{0}, f_{1}, \ldots, f_{k-1}$. Determine $n$ so that

\vspace{+1ex}
$n \equiv f_{i}\ [mod\ (1+(i+1)l!)]$ for $i = 0, 1, \ldots, k-1$,

\vspace{+1ex}
\noindent which is possible, since any two of the numbers $1+(i+1)l!\ (i = 0, 1, \ldots, k-1)$ are relatively prime. For, a prime dividing two of these numbers must also divide the difference $(i_{1} - i_{2})l!$ and therefore, since $i_{1} - i_{2} < l$, must also divide $l!$, which is impossible. The number pair $n, l!$ fulfills our requirement.

\vspace{+1ex}
Since the relation $x = [n]_{p}$ is defined by

\vspace{+1ex}
$x \equiv n\ (mod\ p)\ \&\ x < p$

\vspace{+1ex}
\noindent and is therefore arithmetical, then so also is the relation $P(x_{0}, x_{1}, \ldots, x_{n})$ defined as follows:

\begin{quote}
$P(x_{0}, x_{1}, \dots, x_{n}) \equiv (\exists n, d)\{S([n]_{d+1}, x_{2}, \ldots, x_{n})\ \&\ (\forall k)[k < x_{1} \rightarrow T([n]_{1+d(k+2)}, k,$ $[n]_{1+d(k+1)}, x_{2}, \ldots, x_{n})]\ \&\ x_{0} = [n]_{1+d(x_{1} + 1)}\}$
\end{quote}

\noindent which, according to (17) and Lemma 1, is equivalent to $x_{0} = \phi(x_{1}, \ldots, x_{n})$ (in the sequence $f$ in (17) only its values up to the $(x+1)$th term matter). Thus, Theorem VII(2) is proved.

\begin{quote}
\footnotesize{
Note: G\"{o}del's remark that ``in the sequence $f$ in (17) only its values up to the $(x+1)$th term matter" is significant for the resolution of the PvNP problem. The proof of Lemma \ref{sec:2.3.lem.1}---and consequently of Theorem \ref{sec:2.3.thm.1} that P$\neq$NP if the standard interpretation of PA is sound---depends upon the fact that the equivalence between $f^{(n, d)}$ and $f$ cannot be extended non-terminatingly.
} 
\end{quote}

\section{Appendix C: PA \textit{cannot} admit a set-theoretical model}
\label{sec:5.4.1}

Let $[G(x)]$ denote the PA formula:

\addvspace{+1ex}
$[x=0 \vee \neg(\forall y)\neg(x=y^{\prime})]$

\addvspace{+1ex}
This translates, under every unrelativised interpretation of PA, as:

\addvspace{+1ex}
If $x$ denotes an element in the domain of an unrelativised interpretation of PA, either $x$ is 0, or $x$ is a `successor'.

\addvspace{+1ex}
Further, in every such interpretation of PA, if $G(x)$ denotes the interpretation of $[G(x)]$:

\addvspace{+1ex}
(a)	$G(0)$ is true;

(b)	If $G(x)$ is true, then $G(x^{\prime})$ is true.

\addvspace{+1ex}
Hence, by G\"{o}del's completeness theorem:

\addvspace{+1ex}
(c)	PA proves $[G(0)]$;

(d)	PA proves $[G(x) \rightarrow G(x^{\prime})]$.

\addvspace{+1ex}
\noindent \textit{G\"{o}del's Completeness Theorem}: In any first-order predicate calculus, the theorems are precisely the logically valid well-formed formulas (\textit{i.\ e.\ those that are true in every model of the calculus}).

\addvspace{+1ex}
Further, by Generalisation:

\addvspace{+1ex}
(e)	PA proves $[(\forall x)(G(x) \rightarrow G(x^{\prime}))]$;

\addvspace{+1ex}
\noindent \textit{Generalisation in PA}: $[(\forall x)A]$ follows from $[A]$.

\addvspace{+1ex}
Hence, by Induction:

\addvspace{+1ex}
(f)	$[(\forall x)G(x)]$ is provable in PA.

\addvspace{+1ex}
\noindent \textit{Induction Axiom Schema of PA}: For any formula $[F(x)]$ of PA:

$[F(0) \rightarrow ((\forall x)(F(x) \rightarrow F(x^{\prime})) \rightarrow (\forall x)F(x))]$

\addvspace{+1ex}
In other words, except 0, every element in the domain of any unrelativised interpretation of PA is a `successor'. Further, $x$ can only be a `successor' of a unique element in any such interpretation of PA. 

\subsection{PA and ZF have no common model}
\label{sec:5.4.2}

Now, since Cantor's first limit ordinal, $\omega$, is not the `successor' of any ordinal in the sense required by the PA axioms, and if there are no infinitely descending sequences of ordinals\footnote{cf.\ \cite{Me64}, p261.} in a model---if any---of set-theory, PA and Ordinal Arithmetic\footnote{cf.\ \cite{Me64}, p.187.} cannot have a common model, and so we cannot consistently extend PA to ZF simply by the addition of more axioms.

\subsubsection{\textit{Why} PA has no set-theoretical model}
\label{sec:5.4.3}

We can define the usual order relation `$<$' in PA so that every instance of the Induction Axiom schema, such as, say:

\vspace{+1ex}
(i) $[F(0) \rightarrow ((\forall x)(F(x) \rightarrow F(x')) \rightarrow (\forall x)F(x))]$

\vspace{+1ex}
yields the PA theorem:

\vspace{+1ex}
(ii) $[F(0) \rightarrow ((\forall x) ((\forall y)(y < x \rightarrow F(y)) \rightarrow F(x)) \rightarrow (\forall x)F(x))]$

\vspace{+1ex}
Now, if we interpret PA without relativisation in ZF in the sense indicated by Solomon Feferman\footnote{\cite{Fe92}.} --- i.e., numerals as finite ordinals, $[x']$ as $[x \cup \left \{ x \right \}]$, etc. --- then (ii) always translates in ZF as a theorem:

\vspace{+1ex}
(iii) $[F(0) \rightarrow ((\forall x)((\forall y)(y \in x \rightarrow F(y)) \rightarrow F(x)) \rightarrow (\forall x)F(x))]$

\vspace{+1ex}
However, (i) does not always translate similarly as a ZF-theorem (\textit{which is why PA and ZF can have no common model}), since the following is not necessarily provable in ZF:

\vspace{+1ex}
(iv) $[F(0) \rightarrow ((\forall x)(F(x) \rightarrow F(x \cup \left \{x\right \})) \rightarrow (\forall x)F(x))]$

\vspace{+1ex}
\textit{Example}: Define $[F(x)]$ as `$[x \in \omega]$'. 

A significant point which emerges from the above is that we cannot appeal unrestrictedly to set-theoretical reasoning when studying the foundational framework of PA.

Reason: The language of PA has no constant that interprets in any model of PA as the set \textit{\textit{N}} of all natural numbers.

Moreover, the preceding sections show that the Induction Axiom Schema of PA does not allow us to bypass this constraint by introducing an ``actual" (\textit{or ``completed"}) infinity disguised as an arbitrary constant - usually denoted by $c$ or $\infty$ - into either the language, or a putative model, of PA.

\noindent \tiny{Authors postal address: 32 Agarwal House, D Road, Churchgate, Mumbai - 400 020, Maharashtra, India.\ Email: re@alixcomsi.com, anandb@vsnl.com.}

\end{document}